\documentclass[a4paper,11pt]{article}\usepackage{amssymb}
\usepackage{amsmath}
\usepackage{latexsym,amsmath}
\usepackage{amsmath}

\usepackage{amsfonts}
\usepackage[english]{babel}
\usepackage[dvips]{graphicx}
\usepackage{graphics}
\usepackage{makeidx}
\usepackage[dvips]{color}

\usepackage{amsmath,amssymb}
\def \E {\mathbb{E}}
\def \Et {\mathbb{E}^{\mathcal{F}_t}}
\def \P {\mathcal{P}}
\def \calp {\mathcal{P}}

\def\F {\mathcal{F}}
\def \R {\mathbb{R}}
\def \cald {\mathcal{D}}
\def \dis {\displaystyle}

\def \< {\langle}
\def \> {\rangle}
\newtheorem{theorem}{Theorem}[section]
\newtheorem{lemma}[theorem]{Lemma}
\newtheorem{proposition}[theorem]{Proposition}

\newtheorem{hypothesis}[theorem]{Hypothesis}

\newtheorem{remark}[theorem]{Remark}
\numberwithin{equation}{section}
\def\Dim{\noindent\hbox{{\bf Proof.}$\;\; $}}          
\def\finedim{{\hfill\hbox{\enspace${ \square}$}} \smallskip}    
\def\sqr#1#2{{\vcenter{\vbox{\hrule height .#2pt
     \hbox{\vrule width .#2pt height#1pt \kern#1pt \vrule
     width .#2pt} \hrule height .#2pt}}}}
\def\square{\mathchoice\sqr54\sqr54\sqr{4.1}3\sqr{3.5}3}

\begin{document}
\title{On the existence of optimal controls for SPDEs with boundary-noise and
boundary-control}
\author{Giuseppina Guatteri, \\
Dipartimento di Matematica, \\
Politecnico di Milano, \\
Piazza Leonardo da Vinci 32, \\
20133 Milano. \\ \\
Federica Masiero, \\
Dipartimento di Matematica e Applicazioni, \\
Universit\`{a} di Milano
Bicocca,\\
via R. Cozzi 53 - Edificio U5,\\
20125 Milano.\\ \\
 {\tt e-mail: giuseppina.guatteri@polimi.it, federica.masiero@unimib.it}}
\date{}
\maketitle
{\abstract
We consider a stochastic optimal control problem for an heat equation with
boundary noise and boundary controls. Under suitable assumptions on the coefficients, we
prove existence of optimal controls in strong sense by solving the stochastic hamiltonian system related.  } 
 

\section{Introduction}
In this  paper we are concerned with the existence of optimal control for a stochastic 
optimal control problem related to the following stochastic heat equation, in which
boundary noise and boundary control are allowed:
\begin{equation}\label{eqconcretaintro}
  \left\{
  \begin{array}{l}
  \displaystyle
\frac{ \partial y}{\partial t}(t,\xi)= \frac{ \partial^2
y}{\partial \xi^2}(t,\xi)+b(\xi)u^0(t,\xi)+g(\xi)\dot{W}(t,\xi)  , \qquad t\in [0,T],\;
\xi\in (0,\pi),
\\ \displaystyle
y(0,\xi)=x(\xi),
\\\displaystyle
\frac{ \partial y}{\partial \xi}(t,0)= u^1_s+\dot{\tilde W}_s, \quad
\frac{ \partial y}{\partial \xi}(t,\pi)= u^2_s
\end{array}
\right.
\end{equation}
In the above equation $\tilde W$ is a standard real Wiener process and 
$\dot{W}\left(  \tau,\xi\right)  $ is a space-time white noise on
$\left[  0,T\right]  \times\left[  0,\pi\right]  $; $\tilde W$ and $W$ are
both defined on a complete probability space $(\Omega,\mathcal{F},\P)$
and are independent. By $\{\mathcal{F}_t, \ t \in [0,T] \}$ we will denote the natural filtration
of $(\tilde{W},W)$, completed in the usual way; $u^0$
and $(u^1,u^2)$ are $\mathcal{F}_t$-predictable square integrable processes and represent
respectively the distributed and the boundary control.
Notice that we are able to treat equations where
the control affects all the boundary while the noise only affects one
point at the boundary.

The problem is considered in its {\em strong formulation}, i.e. without changing the reference probability space $(\Omega,\mathcal{F},\P)$.
\noindent The stochastic optimal control problem
consists in minimizing over all admissible controls the following
cost functional:
\begin{multline}\label{costoconcretointro}
J(x,u^0,u^1,u^2)=\E \int_0^T\int_0^{\pi}( \bar l_o(s,\xi,
y(s,\xi))+ \bar g(u^0_s(\xi),u^1_s,u^2_s))\;d\xi\;ds  \\ 
+\E \int_0^{\pi} \bar h(\xi,
y(T,\xi))\;d\xi,
\end{multline}
where $\bar g$ and $ \bar h$ satisfies suitable assumptions specified in section \ref{sez-abstr},
here we only mention that $\bar g$ is allowed to have quadratic growth with respect to the control,
and the control processes are not necessarily bounded.
Equation (\ref{eqconcretaintro}) will be reformulated as a stochastic evolution
equation in $H=L^2((0,\pi))$:
\begin{equation}\label{eqastrattaintro}
\left\{\begin{array}{ll}  dX_t = A X_t \, dt + [(\lambda -A)D + B]u_t \,dt 
+ (\lambda -A)D_1\,d\tilde{W}_t + G\, dW_t & t \in [0,T]\\
X_0 =x, \end{array}\right.
\end{equation} 
where $B$ and $G$ are as usual the multiplication operators related to $b$
and $g$ respectively and $D$ and $D_1$ transform boundary data in elements of the domain of a suitable fractional power of $(\lambda -A)$, so that
both $(\lambda -A)D$ and $(\lambda -A)D_1$ are unbounded operators. Notice that
equation (\ref{eqastrattaintro}) can be considered as the model for a more general class of
state equations, see section \ref{sez-abstr} for more details.

An approach to prove existence of optimal controls is the dynamic programming principle
and the solution, in a sufficiently regular sense, e.g. mild, of the Hamilton Jacobi Bellman (HJB in the following) equation
related. Because of the presence of the boundary noise, the transition semigroup related
to equation (\ref{eqastrattaintro}) does not have sufficient smoothing properties, so the HJB equation associated
cannot be solved in mild sense by a fixed point argument. The HJB equation
is solvable in the sense of viscosity solutions, see e.g. 
\cite{GoRoSw}, and the presence of the noise as a forcing term is necessary in their approach.
Moreover, since in equation (\ref{eqastrattaintro}) the control is not assumed to be
in the image of $G$ nor in the image of $(\lambda -A)D_1$, the HJB equations
cannot be solved by means of backward stochastic differential equations (BSDEs in the following), see the
pioneering paper \cite{PaPe1} and the infinite dimensional extension in
\cite{fute}. When HJB equations can be solved by means of BSDEs, 
boundary noise and boundary control problems for the heat equations
are treated in \cite{DebFuhTes}, in the case of Neumann boundary conditions, and 
the techniques have been extended to the case of Dirichlet boundary conditions in
\cite{Mas}, by using also results in
\cite{FabGol}. 
We also mention that in the dynamic programming approach existence of optimal controls
is proved in the weak sense, since once the HJB
equations is solved, the synthesis of the optimal contros is subject to the solution of the
so called closed loop equation: since it is not clear the regularity of the feedback law, in many cases the closed lopp equation can
be solved only in the weak sense.

\noindent In \cite{GaSo}, by extending finite dimensional techniques,
existence of optimal controls in the case of Hilbert space valued controlled diffusions
is proved in relaxed sense.
In \cite{BuRa} existence of {\em quasi}-optimal controls is proved for a control
problem related to a controlled state equation with distributed control and noise via
the Ekeland principle.  Their setting is infinite dimensional as in the present paper,
but they prove existence of optimal controls not in strong sense
and moreover in the state equation no unbounded terms are allowed.  On the other hand they can bypass convexity assumptions either on the coefficients (still very regular) of the cost functional
or of the control space $U$.

\noindent An other approach to prove existence of optimal controls is the stochastic maximum principle,
see e.g. \cite{HuPeng1}, which provides ncessary conditions for optimality. When these conditions
are also sufficient, existence of optimal controls can be proved by solving
the related forward backward stochastic Hamiltonian system, see e.g. \cite{HuPeng2}.
Both in \cite{HuPeng1} and in \cite{HuPeng2} the setting is finite dimensional
In this paper we generalize this approach to the infinite dimensional setting. The maximum principle,
see \cite{Gua} where the boundary case is treated, provides as usual necessary conditions for the optimal control to be verified. Then, under suitable assumptions -see section 2.2-, one can show that these conditions are indeed sufficient and so the solution to the Hamiltonian system fully characterizes the optimal control.
In our case the Hamiltonian system is the following:
\begin{equation}\label{forwardback_ottimo_intro}
\left\{\begin{array}{ll}  d\bar{X}_t = A \bar{X}_t , dt + [E + B]\gamma([E+ B]^*\bar{Y}_t)\,dt+ (\lambda -A)D_1\,d\tilde{W}_t + G(t,\bar{X}_t)\, dW_t &\\ \\
 - d{\bar{Y}}_t = A^* {\bar{Y}}_t \, dt   + l^0_x(t,\bar{X}_t) \, dt  -\bar{Z}_t \, dW_t  -\tilde{Z}_t \, d\tilde{W}_t , \qquad t \in [0,T] & \\ \\
\bar{X}_0 =x, \ \bar{Y}_T = - h_x(\bar{X}_T) ,\end{array}\right.
\end{equation}
where $H(t,x,u,y):= - l(t,x,u) + \langle  [E + B]^*y, u\rangle,$ is the hamiltonian function, and $\gamma: H \to U$ is such that 
$H(t,x,\gamma([E+B]^*y),y) = \inf _{u \in U } H(t,x,u,y)$.
Because of the infinite dimensional setting and of the presence of unbounded operators, the result obtained
in the solution of this infinite dimensional forward backward system are of independent interest.

Indeed the solution of fully coupled forward backward systems is a
difficult topic already in the finite dimensional case, see \cite{Anto} and again \cite{MaYong}
for examples of finite dimensional FBSDEs where there is no hope to get existence of a solution.

Among the large literature in finite dimensions, see e.g. the book of \cite{MaYong},
we can distinguish two main approaches. The first approach, known as \textit{four-step schem}, relies on
the connections between SDEs with deterministic coefficients and non-linear PDEs, see the
pioneering paper \cite{MaProYong}. 
Since in infinite dimensions on the solution of the related PDE less apriori estimates are known, 
this approach seems to be not suitable for an infinite dimensional extensions: in \cite{Gua}
local existence for an infinite dimensional FBSDE is proved, mainly adequating the finite dimensional techniques introduced
in \cite{Delarue}, but global existence is not achieved.

The second approach applies
under monotonicity assumptions: different types of
conditions have been investigated in this framework and we refer to Hu and Peng \cite{HuPeng3},
Peng and Wu \cite{peng wu 1999}, Yong \cite{yong 1997}
and to Pardoux and Tang \cite{pardoux tang 1999}.
\vspace{5pt}
\\
%

In the present paper, we solve FBSDE (\ref{forwardback_ottimo_intro}) by adapting the
\emph{bridge method} introduced in \cite{HuPeng3} to the infinite dimensional framework:  new difficulties arises because of the presence of the
unbounded operators, and just because both the forward and the backward stochastic equations are infinite dimensional and an unbounded operator is applied to backward unknown $Y$ in the forward equation so that one has to prove some extra regularity for $Y$ in order to give meaning to the system in the space $H$.
The regularity of the adjoint unknown is a typical task when one wants to prove maximum principle in infinite dimension, see \cite{HuPeng1} and \cite{Gua},  in this case new difficulty arise since the backward equation is coupled with the forward and the whole system has to be considered.
The linear auxiliary FBSDE we study to apply then the \emph{bridge method} is
\begin{equation}\label{forwardback_linear_bridge_intro}
\left\{\begin{array}{ll}  d\bar{X}_t = A \bar{X}_t , dt - [E + B][E+ B]^*\bar{Y}_t\,dt+ b_0(t)\,dt + (\lambda -A)D_1\,d\tilde{W}_t + G\, dW_t &\\ \\
 - d{\bar{Y}}_t = A^* {\bar{Y}}_t \, dt   +\bar{X}_t \, dt + h_0(t) \, dt  -\bar{Z}_t \, dW_t  -\tilde{Z}_t \, d\tilde{W}_t , \qquad t \in [0,T] & \\ \\
\bar{X}_0 =x, \, -\bar{Y}_T = \bar{X}_T + g_0, \end{array}\right.
\end{equation}

Unlike in \cite{HuPeng3}, this linear auxiliary FBSDE is not immediately solvable. We
notice that such system is the hamiltonian system associated to of an affine quadratic
optimal control problem with state equation
\begin{equation}\label{forward_lin.intro}
\left
\{\begin{array}{ll}  dX_t = A X_t \, dt + [E+ B]u_t \,dt +  b_0(t) \, dt + (\lambda -A)D_1\,d\tilde{W}_t + G\, dW_t & t \in [0,T]\\
X_0 =x ,\end{array}
\right.
\end{equation} 
and cost functional 
\begin{equation}\label{funzionalecosto.intro}
J(x,u)= \frac{1}{2} \E \int_0^T (|X_t + h_0(t)|^2 + |u_t|^2) \, dt +  \frac{1}{2} \E |X_T+g_0|^2
\end{equation}
where $b_0$ ad $h_0$ are suitable stochastic processes.
Therefore we introduce the Riccati equation (deterministic) corresponding to the linear terms and a backward stochastic differential equation to deal with the affine terms, see section 3.2, in order to get a solution to system \eqref{forwardback_linear_bridge_intro}. Again, because of the infinite dimensional setting and of the presence of unbopunded operators, the solution and the regularity
of this auxiliary backward stochastic differential equation is of independent interest.

\noindent Once we prove that  system \eqref{forwardback_linear_bridge_intro} has a unique solution, for every suitable $b_0$ and $h_0$ we can start to "build" the {\em bridge} to get a solution to our original system \eqref{forwardback_ottimo_intro} and then eventually solve our control problem.

The paper is organized as follows; in section 2 we state the notations and the problem,
we collect results on the stochastic maximum principle in the boundary case and we also prove suffiecient 
contitions for otpimality, finally we stae our main result on the existence of optimal controls; in
section 3 we prove existence and uniqueness of a mild solution for the stochastic hamiltonian
system by applying the bridge method to this setting and we conclude by proving the existence of optimal controls.

\section{Preliminaries and statement of the problem}
\subsection{Notation}
Given a Banach space $X$,
the norm of its elements $x$ will be denoted
by $|x |_X$, or even by $|x|$ when no confusion is possible. If $V$ is
another Banach space, $L(X,V)$ denotes the space of bounded linear
operators from $X$ to $V$, endowed with the usual operator norm.
Finally we say that a mapping $F : X \to V$ belongs to the class $\mathcal{G}^1(X;V)$
if it is continuous, G\^ateaux differentiable on $X$,
 and $\nabla F: X\to L(X,V)$ is strongly continuous.
The letters $\Xi$, $H$, $K$ and $U$ will always be used to denote Hilbert spaces.
The scalar product is denoted $\langle \cdot, \cdot \rangle$, equipped with a
subscript to specify the space, if necessary. All the Hilbert
spaces are assumed to be real and separable; $L_2(\Xi,H)$ is the space of Hilbert-Schmidt operators from $\Xi $ to $H$, respectively.

\vspace{5pt}
 Given an arbitrary but fixed time horizon  $T$, we consider all stochastic processes
 as defined on subsets of the time interval $[0,T]$.
 Let $Q \in L(K)$ be a symmetric non-negative operator, not necessarily trace class and $\tilde{W} = (\tilde{W}_t)_{t \in [0,T]}$  be a $Q$-Wiener process with
values in $K$, defined on a
complete
probability space $(\Omega, \mathcal{F}, \mathbb{P})$ and $W = (W_t)_{t \in [0,T]}$ be a cylindrical Wiener process with values in $\Xi$, defined on the same probability space and independent of $\tilde{W}$.
By $\{\mathcal{F}_t, \ t \in [0,T] \}$ we will denote the natural filtration
of $(\tilde{W},W)$, augmented with the family $\mathcal{N}$ of
$\mathbb{P}$- null sets of $\mathcal{F}$, see for instance \cite{DPZ1} for its definition. Obviously, the filtration
$(\mathcal{F}_t)$ satisfies the usual conditions of right-continuity and completeness. All the concepts
of measurability for stochastic processes will refer to this filtration.
By $\mathcal{P}$ we denote the predictable $\sigma$-algebra on
$\Omega \times [0,T]$ and by $\mathcal{B}(\Lambda)$ the Borel
$\sigma$-algebra of any topological space $\Lambda$.

Next we define two classes of stochastic processes with values in
a Hilbert space $V$.
\begin{itemize}
\item $L^2_\calp (\Omega\times [0,T];V)$ denotes the space of
equivalence classes of processes $Y \in L^2 (\Omega\times
[0,T];V)$ admitting a predictable version. It is endowed with the norm
\[ |Y|= \Big(\E \int_0^T |Y_s|^2 \, ds\Big)^{1/2}. \]
 \item $C_\calp([t,T];L^p(\Omega;S))$, $p\in [1,+\infty]$, $t \in [0,T]$, denotes the space of $S$-valued processes $Y$ such that
 $Y : [t,T] \to L^p(\Omega,S)$ is continuous
 and $Y$ has a predictable modification, endowed with the norm:
\begin{equation*}
|Y|^p_{C_\calp([t,T];L^p(\Omega;S))}=\sup_{s \in [t,T]}\E
|Y_s|^p_S
\end{equation*}
Elements of $C_\calp([t,T];L^p(\Omega;S))$ are identified up to
modification.
\item For a given $p \geq 2$, $L^p_{\mathcal{P}}(\Omega;C([0,T];V))$
     denotes the space of
    predictable processes $Y$ with continuous paths in $V$, such
    that the norm
    \[  \|Y\|_p  = (\E \sup _{s \in [0,T]} |Y_s|^p)^{1/p}\]
    is finite. The elements of $L^p_{\mathcal{P}}(\Omega;C([0,T];V))$
    are identified up to indistinguishability.
\end{itemize}
Given an element $\Phi$ of  $L^2_\P (\Omega\times
[0,T];L_2(\Xi,V))$ or of $L^2_\P (\Omega\times
[0,T];L_2(K,V))$, the It\^o stochastic integrals $\int_0^t
\Phi(s) \,dW(s)$ and  $\int_0^t \Phi(s) \,d\tilde{W}(s)$, $t \in
[0,T]$, are  $V$-valued martingales belonging to
$L^2_{\mathcal{P}}(\Omega;C([0,T];V))$. The previous definitions
have obvious extensions to processes defined on subintervals of
$[0,T]$ or defined on the entire positive real line $\R^+$.
\subsection{Optimal control problem and state equation}
\label{sez-abstr}
Let $H$ be a separable real Hilbert space, and $U$ a separable Hilbert, called the space of controls.
We assume $U$ a convex set and we set the space $L^2_\P (\Omega\times [0,T];U)$ the space of admissible controls, and we denote it by $\mathcal{U}$.

We make the following assumptions that we denote by ${\bf(A)}$:
\begin{enumerate}
\item[{\bf (A.1)}] $A: D(A)\subset H \to H$ is a 
linear, unbounded operator that generate a $C_0$-semigroup $\{
e^{tA} \}_{t \geq 0}$ that is also analytic and 
such that $|e^{tA}|_{L(H,H)} \le  e ^{ \omega t}$, $ t \geq 0$
 for some  $\omega \in \R$. This means in particular that every $\lambda > \omega $ belongs to the resolvent set of $A$. \\
\item[{\bf(A.2)}] $B \in L(U;H)$
and 
 $G \in L(\Xi,H)$ and there exist constants  $\Delta > 0$ and $\gamma \in
[0,1/2[$  such that
$$
|e^{sA}G|_{L_2(\Xi,H)}  \leq \frac{\Delta }{(1\wedge
s)^\gamma}$$
 for every $s \in \mathbb{R}^+$.
 \item[{\bf (A.3}] $D$ is a continuous linear operator  $D: U \to D((\lambda -A )^\alpha)$ for some $\frac{1}{2} < \alpha <1$ and $\lambda > \omega$, see for instance \cite{Lunardi} or \cite{Pazy} for the definition of the fractional power of the operator $A$. 
\item[{\bf (A.4)}] $D_1$ is a linear operator $D_1: K\to H$ and there is a constant $\frac{1}{2}< \beta < 1 $ such that the following holds:
\[  |e^{tA}(\lambda-A)D_1|_{L_2(K,H)} \leq \frac{C}{t^{1-\beta}} \] 
 \end{enumerate}
 for some $\lambda >0$.
\begin{remark}
 Notice that $D_1$ and $D$ can have the same structure, indeed if  $D_1$ takes values in 
$D((\lambda -A )^\beta)$ and $K$ is finite dimensional then ${\bf (A.4)}$ holds.  On the over hand, by the analiticity of $A$,  also for $D$ a similar estimate to the one for
$D_1$ may follow.
\end{remark}
\vspace{2em}

We introduce the following  class  of control problems, where the {\em state equation} is 
\begin{equation}\label{forward}
\left\{\begin{array}{ll}  dX_t = A X_t \, dt + [(\lambda -A)D + B]u_t \,dt + (\lambda -A)D_1\,d\tilde{W}_t + G\, dW_t & t \in [0,T]\\
X_0 =x \end{array}\right.
\end{equation} 
From now on we will denote for semplicity $(\lambda -A)D :=E$

We will seek for a {\em mild} solution to this equation, in the sense of \cite{DPZ1}, that is  a $(\F_t)$- predictable process $X_t, \ t \in [0,T]$  with continuous path in $H$ such that $\P$- a.s.
\begin{multline}\label{stato_mild}
X_t = e^{tA} x + \int _0^t e^{(t-s)A}[E +B]u_s \,ds  + \int _0^t e^{(t-s)A}(\lambda -A)D_1 \, d\tilde{W}_s   +\int_0^t e^{(t-s)A} G\, dW_s, \ t \in [0,T]
\end{multline}
The {\em cost functional}, that depends on the initial state $x$ and the control $u \in \mathcal{U}$, to minimize is:
\begin{equation}\label{funzionale}
J(x,u)= \E \int_0^T l(t,X_t,u_t) \, dt + \E h(X_T)
\end{equation}
where $l$ and $h$ verify {\bf(B)}:
\begin{enumerate}
\item[{\bf(B.1)}] $l$ is measurable and for all $ t\in [0,T]$ and all $u \in U$, $l(t,\cdot,u) \in \mathcal{G}^1(H;\R)$ and for all $ t\in [0,T]$ and all $x \in H$, $l(t,x,\cdot) \in \mathcal{G}^1(U;\R)$ and there is a constant $\Delta >0$ such that:
\begin{equation}
|l_x(t,x,u)|+ |l_u(t,x,u)| \leq \Delta ( 1+ |x|_H + |u|_U) 
\end{equation}
for all $t \in [0,T]$, $x \in H$ and $u \in U$.
\item[{\bf(B.2)}] the map $h$ is continuous and convex, moreover $h \in \mathcal{G}^1(H;\R)$ and there is a constant $\Delta >0$ such that:
\begin{equation}
|h_x(x)| \leq \Delta ( 1+ |x|_H ) 
\end{equation}
for all $x \in H$.
Moreover  for some constant $c_1 >0$
\begin{equation}\label{dissipativita_datofin}
\langle h_x(x_1) - h_x(x_2), x_1-x_2  \rangle _H  \leq -c_1 |x_1-x_2|^2, \text{ for any } x_1, x_2 \in H 
\end{equation}
\item[{\bf(B.3)}]  the map $l$ can be decomposed as $l(t,x,u)=l^0(t,x) + g(u)$ , where $l^0$ and $g$ are two convex functions. 
Moreover  for some constant $c_1 >0$
\begin{equation}\label{dissipativita_costo}
\langle l^0_x(t,x_1) - l^0_x(t,x_2), x_1-x_2  \rangle _H  \geq c_1 |x_1-x_2|^2, \text{ for any } x_1, x_2 \in H, t \in [0,T] 
\end{equation}
\item[{\bf(B.4)}]  for any $t \in [0,T], x \in H,  y \in D(E^* ) $,  we define  $$H(t,x,u,y):= - l(t,x,u) + \langle  [E + B]^*y, u\rangle, $$  and assume that there exists a function $\gamma: H \to U$ such that 
\begin{equation}\label{defdigamma}
H(t,x,\gamma([E+B]^*y),y) = \inf _{u \in U } H(t,x,u,y).
\end{equation}
We assume moreover that there exist positive constants $c_1$ and $\Delta$:
\begin{equation}\label{dissipativita_gamma}
\langle \gamma (y_1) -\gamma (y_2), y_1-y_2  \rangle _H  \leq -c_1 |y_1-y_2|^2, \text{ for any } y_1, y_2 \in H 
\end{equation}
\begin{equation}\label{Lip_gamma}
|\gamma (y_1) -\gamma (y_2)|_H \leq  \Delta |y_1-y_2|, \text{ for any } y_1, y_2 \in H 
\end{equation}
\end{enumerate}

\subsection{Heat Equation with Neumann Boundary conditions}
\label{sez-heateq}

In this section we present a concrete stochastic control problem
that we will be able to treat and we show how this model fits the ``abstract''
setting of section \ref{sez-abstr}.
We consider an heat equation on the interval $(0,\pi)$
with boundary noise and boundary control, and we focus our attention
on the case where the control affects all the boundary, and the noise affects only
one point at the boundary.
\begin{equation}\label{eqconcreta}
  \left\{
  \begin{array}{l}
  \displaystyle
\frac{ \partial y}{\partial t}(t,\xi)= \frac{ \partial^2
y}{\partial \xi^2}(t,\xi)+b(\xi)u^0(t,\xi)+g(\xi)\dot{W}(t,\xi)  , \qquad t\in [0,T],\;
\xi\in (0,\pi),
\\ \displaystyle
y(0,\xi)=x(\xi),
\\\displaystyle
\frac{ \partial y}{\partial \xi}(t,0)= u^1_t+\dot{\tilde W}_t, \quad
\frac{ \partial y}{\partial \xi}(t,\pi)= u^2_t
\end{array}
\right.
\end{equation}
In the above equation $\tilde W$ is a standard real Wiener process and 
$\dot{W}\left(  \tau,\xi\right)  $ is a space-time white noise on
$\left[  0,T\right]  \times\left[  0,\pi\right]  $; $\tilde W$ and $W$ are independent.
We will give sense to the
notion of solution in the following.

We reformulate equation (\ref{eqconcreta}) as a stochastic evolution
equation in $H=L^2(0,\pi)$.
$A$ stands for the Laplace operator with
homogeneous Neumann boundary conditions,
which is the generator of an analytic semigroup in $H$:
$$
\cald(A)=\left\lbrace y\in H^2(0,\pi):\frac{ \partial y}{\partial \xi}(0)=
\frac{ \partial y}{\partial \xi}(\pi)=0\right\rbrace, \qquad 
Ay=\frac{ \partial^2 y}{\partial \xi^2} \text{ for }y\in \cald(A).
$$
The control process $u\in L_\P^2(\Omega\times[0,T],U)$ where
$U=L^2(0,\pi)\times \R^2$ and $u=\left(\begin{array}{l}
u^0\\u^1\\u^2
\end{array}\right)$
We fix $\lambda>0$ and define
\[
 b^1(\xi)=-\dfrac{\cosh(\sqrt{\lambda}(\pi-\xi))}{\sqrt{\lambda}\sinh(\sqrt{\lambda}\xi)},\qquad
 b^2(\xi)=\dfrac{\cosh(\sqrt{\lambda}\xi)}{\sqrt{\lambda}\sinh(\sqrt{\lambda}\xi)}
\]
and note that they solve the Neumann problems
\begin{equation*}
  \left\{
  \begin{array}{l}
  \displaystyle
\frac{ \partial^2
b^i}{\partial \xi^2}(\xi)=\lambda b^i(\xi)  , \qquad \xi\in (0,\pi),\;i=1,2,
\\ \displaystyle
\frac{ \partial b^1}{\partial \xi}(0)=1,
\quad \frac{ \partial b^1}{\partial \xi}(\pi)=0\\ \dis
\frac{ \partial b^2}{\partial \xi}(0)=0,
\quad \frac{ \partial b^1}{\partial \xi}(\pi)=1.
\end{array}
\right.
\end{equation*}
So $b^i\in \cald(\lambda-A)^\alpha=H^{2\alpha}$, for $1/2<\alpha<3/4$.

Equation \ref{eqconcreta} can now be reformulated as:
\begin{equation}\label{eqastratta}
\left\{\begin{array}{ll}  dX_t = A X_t \, dt + [(\lambda -A)D + B]u_t \,dt + (\lambda -A)D_1\,d\tilde{W}_t + G\, dW_t & t \in [0,T]\\
X_0 =x, \end{array}\right.
\end{equation} 
where, for $u\in U$ and $h\in H$, $Du=(0,b^1(\cdot)u^1(\cdot), b^2(\cdot)u^2(\cdot))$, $D_1=(0,b^1(\cdot)u^1(\cdot),0)$, $B=(b(\cdot),0,0)$, $Gh=g(\cdot)h(\cdot)$.
With the notations of section \ref{sez-abstr}, $K=\R$ and $\Xi=H$.

\noindent Equation (\ref{eqastratta}) is still formal, since $(\lambda -A)D$ and
$(\lambda -A)D_1$ do not take their values in $H$, the precise meaning of equation (\ref{eqastratta})
is given by its mild formulation.
An $H$-valued predictable process $X$  is called a mild solution 
to equation (\ref{eqastratta}) on $[0,T]$ if 
\[
 \mathbb P\int_0^T\vert X_r\vert^2 dr <+\infty
\]
and, for every $0< t<T$, $X$ satisfies the integral equation
\[
 X_t=e^{tA} x+\int_0^t e^{(t-r)A} [\lambda -A)D+B]u_r dr 
+\int_0^t e^{(t-r)A}(\lambda -A)D_1\,d\tilde{W}_r +\int_0^t e^{(t-r)A} G dW_r.
\]
Since $b^i\in \cald(\lambda-A)^\alpha=H^{2\alpha}$, for $1/2<\alpha<3/4$, and by the analyticity of the semigroup
$e^{tA},\;t\geq 0$, the integral
$\int_0^t e^{(t-r)A}(\lambda -A)Du_r dr$ and the stochastic integral
$\int_0^t e^{(t-r)A}(\lambda -A)D_1\,d\tilde{W}_r$ are well defined, see also \cite{DebFuhTes}.

Notice that equation (\ref{eqastratta}) does not satisfy any structure condition
suitable to treat the related  stochastic optimal control problem 
using backward stochastic differential equations, as in \cite{DebFuhTes}  and \cite{Mas}, where the case of an heat
equation with Dirichlet boundary-control and boundary noise is considered.
Notice that in the present example, differently from   \cite{DebFuhTes} and \cite{Mas}, the control affects the system in $0$ and $\pi$ and the noise
acts only at $0$, so that $\operatorname{Im} (D)\nsubseteq \operatorname{Im} (D_1)$. 

The optimal control problem we wish to treat in this paper
consists in minimizing the following finite horizon cost
 \begin{equation}\label{costoconcreto}
J(x,u^0,u^1,u^2)=\E \int_0^T\int_0^{\pi} \bar l(s,\xi,
y(s,\xi),u^0_s(\xi),u^1_s,u^2_s)\;d\xi\;ds +\E \int_0^{\pi} \bar h(\xi,
y(T,\xi))\;d\xi,
\end{equation}
over all admissible controls. The cost functional \eqref{costoconcreto}
can be written in an abstract way as in (\ref{funzionale}) by setting,
for $s\in[0,T], x\in H, u\in U$
\[
 l(s,x,u)=\int_0^{\pi}l(s,\xi, x(\xi), u^0_s(\xi),u^1_s,u^2_s)\quad h(x)=\int_0^{\pi}\bar h(\xi,x(\xi)).
\]
We consider costs such that $\bar l(s,\xi,y,u^0,u^1,u^2)=\bar l^0(s,\xi,y)+\bar g(\xi,u^0,u^1,u^2)$
so that $l$ can be decomposed as in (\textbf{B.3}).
From $\bar l^0$ and $\bar g$ we define $l^0$ and $g$
as we have defined $l$:
\[
 l^0(s,x)=\int_0^{\pi}\bar l^0(s,\xi, x(\xi))d\xi\quad g(u)=\int_0^{\pi}\bar g(\xi,u^0(\xi),u^1,u^2)d\xi.
\]
We make suitable assumptions on $\bar l^0, \, \bar g,\, \bar h$ such that
$l^0$, $g$ and $h$ satisfy assumptions \textbf{B1}-\textbf{B3}.
\begin{hypothesis}\label{ip_costo_concreto}
We assume that:
\begin{itemize}
 \item[1)]the map $\bar h:[0,\pi]\times\R\rightarrow\R$, is measurable, for a.a. $\xi\in[0,\pi]$
$\bar h(\xi,\cdot):\R\rightarrow\R$  is continuous, convex and differentiable
and there exists $\Lambda \in L^{\infty}([0,\pi])$ such that
\[
\vert h_x(\xi, x)\vert\leq \Lambda(\xi)(1+\vert x\vert).
\]
Moreover we assume that for a.a. $\xi\in[0,\pi]$
$\bar h(\xi,\cdot):\R\rightarrow\R$
is dissipative, namely, for every $x_1,x_2\in \R$
\[
 (\bar h(\xi,x_1)-\bar h(\xi,x_2))(x_1-x_2)\leq -c_1(x_1-x_2)^2;
\]
for some positive constant $c_1$.
\item[2)]the map $\bar l^0:[0,T]\times[0,\pi]\times\R\rightarrow\R$ is measurable and for a.a.
$t\in[0.T]$ and $\xi\in[0,\pi]$, $\bar l^0(t,\xi,\cdot)$ is continuous, convex and differentiable,
and there exists $\Lambda \in L^\infty([0,\pi])$ such that $\forall \xi\in[0,\pi]$ and $\forall x\in\R$
\[
\vert \bar l ^0_x(t,\xi, x)\vert\leq \Lambda(\xi)(1+\vert x\vert).
\]
Moreover we assume that for a.a. $t\in[0,\pi]$ and $\xi\in[0,\pi]$
$\bar l^0(t,\xi,\cdot):\R\rightarrow\R$
is dissipative, namely, for every $x_1,x_2\in \R$
\[
 (\bar l^0(t,\xi,x_1)-\bar l^0(t,\xi,x_2))(x_1-x_2)\geq - c_1(x_1-x_2)^2;
\]
for some positive constant $c_1$.
\item[3)] the map $\bar g:[0,\pi]\times\R\times\R\times\R\rightarrow\R$
is measurable and for a.a. $\xi\in[0,\pi]$, $\bar g(\xi,\cdot,\cdot,\cdot):\R^3\rightarrow\R$
is continuous, convex and differentiable and there exists $\Lambda \in L^\infty([0,\pi])$ such that
\[
\vert \bar g_{u^0}(\xi,u^0,u^1,u^2)\vert\leq \Lambda(\xi)(1+\vert u^0\vert)
\]
and a constant $c>0$ such that
\[
\vert \bar g_{u^i}(\xi,u^0,u^1,u^2)\vert\leq c(1+\vert u^i\vert), \qquad i=1,2.
\]
\end{itemize}
\end{hypothesis}

\section{Main results}
In this section we come back to the abstract formulation of the problem, introducing the scheme we follow to find the optimal control: first we prove the maximum principle, then we prove that under our assumptions the condition is also sufficient and in the end we introduce the Hamiltonian system to be solved.
\subsection{Maximum principle}
Let us assume that there exists an optimal control $\bar{u} \in \mathcal{U}$,
under hypotheses stated previously we have that there exists a unique mild solution ${\bar{X}}$ to \eqref{forward} corresponding to $\bar{u}$, see for instance \cite{DPZ1}.
 So $(\bar {u},\bar{X})$ is an {\em optimal pair} for the control problem described by \eqref{forward} and  \eqref{funzionale}.
We introduce the following {\em forward-backward} system, composed by the state equation corresponding to the optimal control $\bar{u}$ and its adjoint equation:
\begin{equation}\label{sistema}
\left\{\begin{array}{ll}  d\bar{X}_t = A \bar{X}_t , dt + [E + B]\bar{u}_t\,dt+ (\lambda -A)D_1\,d\tilde{W}_t + G\, dW_t &\\ \\
 - d{\bar{Y}}_t = A^T {\bar{Y}}_t \, dt   + l^0_x(t,\bar{X}_t) \, dt  -\bar{Z}_t \, dW_t  -\tilde{Z}_t \, d \tilde{W}_t, \qquad t \in [0,T] & \\ \\
\bar{X}_0 =x, \ \bar{Y}_T = - h_x(\bar{X}_T) \end{array}\right.
\end{equation}
Once the forward equation is solved, the adjoint equation is a backward equation depending on the parameter $\bar{X}$.  The existence and uniqueness of a {\em mild} solution  $(\bar{Y},(\bar{Z},\tilde{Z}))  \in L^2_\P(\Omega;C([0,T];H)) \times  L^2_\P ([0,T]\times \Omega; L_2(\Xi \times K,H) )$  for such equation was firstly proved in see \cite{HuPeng}.
We collect  the mentioned results in this proposition:
\begin{proposition}
Assume {\bf (A)} and {\bf (B)}. System \eqref{sistema} has a unique mild solution $(\bar{X},\bar{Y},\bar{Z})$. Moreover  \begin{equation}
\sup_{t \in [0,T[} \E(T-t)^{2(1-\alpha)}\|\bar{Y_t}\|^2_{D(E^*)} < + \infty 
\end{equation}\end{proposition}
\Dim The regularity result  can be proved as in proposition 3.1 of \cite{Gua}.
\finedim
\begin{theorem} \label{massimo_convesso} Assume {\bf (A)} and {\bf (B)}.
Let $(\bar {u},\bar{X})$ be an optimal pair for the problem  \eqref{forward} and  \eqref{funzionale}. Then there 
exists a unique pair $(\bar{Y},(\bar{Z},\tilde{Z})) \in L^2_\P(\Omega;C([0,T];H)) \times  L^2_\P ([0,T]\times \Omega; L_2(\Xi \times K,H) )$ solution of equation \eqref{sistema} such that: 
\begin{equation}\label{hamiltoniana}
\langle H_u(t, \bar{X}_t,\bar{u}_t,Y_t), v -\bar{u}_t \rangle \leq 0, \qquad \forall v \in U,  \ a.e. \ t \in [0,T], \ \mathbb{P}-a.s.
\end{equation}
where 
\begin{equation*}
H(t,x,u,p):=
   \langle (E+B)^*p,u\rangle_H -l(t,x,u), \qquad (t,x,u,p) \in [0,T] \times H \times U \times D(E^*), \ \lambda > \omega
\end{equation*}
\end{theorem}
\Dim The result follows from theorem 4.6 of \cite{Gua} taking $F_x$  and $G_x$ equal to zero; the presence of the bounded operator $B$ does not introduce any new difficulty. The proof follows exactly in the same way. \finedim

\subsection{Sufficient condition for optimality}
Now we present the following sufficient condition of optimality.
Let us consider the forward-backward system \eqref{sistema}: for any admissible control $\bar{v} \in \mathcal{U}$ there exists a solution $(\bar{X},\bar{Y},(\bar{Z},\tilde{Z}))$ we say then 
that $(\bar {v},\bar{X},\bar{Y},(\bar{Z},\tilde{Z}))$ is an {\em admissible 4-tuple}. 
\begin{theorem}\label{condsuff}
Assume {\bf (A)} and {\bf (B)}. Let $(\bar {u},\bar{X},\bar{Y},(\bar{Z},\tilde{Z}))$ an admissible 4-tuple. If
\begin{equation}\label{condizione}
\langle H_u(t, \bar{X}_t,\bar{u}_t,\bar{Y}_t), v -\bar{u}_t \rangle \leq 0, \qquad \forall v \in U,  \ a.e. \ t \in [0,T], \ \mathbb{P}-a.s.
\end{equation}
then  $(\bar {u},\bar{X})$ is optimal for problem \eqref{forward} and  \eqref{funzionale}.
\end{theorem}
\Dim
Let $\bar{v} \in \mathcal{U}$ hence $\bar{u} +\lambda (\bar{v}- \bar{u}) \in \mathcal{U}$, for all $\lambda \in [0,1]$. Being the state equation affine, we have that $\bar{X}^{\bar{u} +\lambda (\bar{v}- \bar{u})}= \bar{X}+ \lambda \tilde{X}^{\bar{v}- \bar{u}}$,
where $\tilde{X}^{\bar v-\bar u}$ solves
the following equation
\begin{equation*}
 \left\lbrace
\begin{array}[c]{l}
 d\tilde{X}^{\bar v-\bar u}_t=A\tilde{X}^{\bar v-\bar u}_t\,dt+(E+B)(\bar{v}_t-\bar{u}_t)\, dt\\
\tilde{X}^{\bar v-\bar u}_0=0,
\end{array}\right.
\end{equation*}
that is, in mild form,
\begin{equation}\label{prima_approx}
\tilde{X}^{v-u}_t=  \int_0^t e^{(t-s)A} (E+B) (\bar{v}_s- \bar{u}_s)\, ds
\end{equation}
Therefore by the convexity assumption of $l_0,g,h$ we end up with:
\begin{align*}
J(x,\bar{u})-J(x,\bar{u}+ \lambda (\bar{v}-\bar{u}))= &\E \int_0^T [l(t,\bar{X}_t,\bar{u}_t)-l(t,\bar{X}_t + \lambda \tilde{X}^{\bar{v}- \bar{u}}_t, \bar{u} _t+\lambda (\bar{v}_t- \bar{u}_t))] \,dt  \\ + &  \ \E [h(\bar{X}_T)- h(\bar{X}_t + \lambda \tilde{X}^{\bar{v}- \bar{u}}_T)] \\
 \leq 
&-\E \int_0^T\lambda \langle l_x^0(t,\bar{X}_t), \tilde {X}^ {\bar{v}-\bar{u} }_t\rangle \,dt - \E \int_0^T \lambda \langle g_u(\bar{u}_t), \bar{v}_t-\bar{u}_t \rangle \, dt
\\
& -\E \lambda \langle  h_x(\bar{X}_T), \tilde{ X}_T ^{\bar{v}-\bar{u}}\rangle .
\end{align*}
Now following the usual approximation strategy we multiply both equations for $\tilde{X}^{\bar{v}-\bar{u}}$
and $\bar{Y}$ by $n (n - A)^{-1}=nR(n,A)$ for $ n > \lambda$, so that the two processes
$\tilde{X}^{{\bar{v}-\bar{u}},n}:= nR(n,A) \tilde{X}^{\bar{v}-\bar{u}}$ and $\bar{Y}^n:=  nR(n,A)\bar{Y}$  both admit an It\^o differential:
\begin{align*}
d\langle \tilde{X}^{{\bar{v}-\bar{u}},n}_t ,\bar{Y}^n_t \rangle= \langle  \bar{Y}^n_t, nR(n,A) [E + B] (\bar{v}_t-\bar{u}_t)\rangle \,dt + \langle  nR(n,A)l_x^0(t,\bar{X}_t), \tilde {X}^ {\bar{u}-\bar{v},n }_t\rangle \,dt
\end{align*}
Observing that $D (E^*) \equiv  D((\lambda - A^* )^{1-\alpha})$, we can let $n$ tend  $\infty$ and we get that:
\begin{align*}
 -\E  \langle  h_x(\bar{X}_T), \tilde{ X}_T ^{\bar{v}-\bar{u}}\rangle -\E \int_0^T\langle l_x^0(t,\bar{X}_t), \tilde {X}^ {\bar{v}-\bar{u} }_t\rangle \,dt =
 \E \int_0^T\langle (v_t-u_t), [E + B]^* \bar{Y}_t\rangle \,dt 
\end{align*}
Notice that
\[
 \E \int_0^T\langle nR(n,A)[E + B](v_t-u_t),  \bar{Y}^n_t\rangle \,dt
\]
makes sense since $\bar{Y}^n_t=nR(n,A)\bar{Y}_t$ and $\bar{Y}_t\in D(E)=D((\lambda-A)^{1-\alpha})$, so also
 $\bar{Y}^n_t\in D(E)$ and also $nR(n,A^*)\bar{Y}^n_t\in D(E)$. So
\[
 \E \int_0^T\langle nR(n,A)[E + B](v_t-u_t),  \bar{Y}^n_t\rangle \,dt=
\E \int_0^T\langle(v_t-u_t), [E + B]^* nR(n,A^*)\bar{Y}^n_t\rangle \,dt
\]
is well defined.
Now, thanks to \ref{hamiltoniana}, taking $\lambda =1$ we get:
\begin{equation*}
J(x,\bar{u}) \leq J(x,\bar{v}),\qquad\qquad \text{ for all } \bar{v} \in \mathcal{U}. 
\end{equation*}
\finedim
\subsection{Hamiltonian System}
Let us introduce the Hamiltonian system associated to our control problem
\begin{equation}\label{forwardback_ottimo}
\left\{\begin{array}{ll}  d\bar{X}_t = A \bar{X}_t , dt + [E + B]\gamma([E+ B]^*\bar{Y}_t)\,dt+ (\lambda -A)D_1\,d\tilde{W}_t + G(t,\bar{X}_t)\, dW_t &\\ \\
 - d{\bar{Y}}_t = A^T {\bar{Y}}_t \, dt   + l^0_x(t,\bar{X}_t) \, dt  -\bar{Z}_t \, dW_t  -\tilde{Z}_t \, d\tilde{W}_t , \qquad t \in [0,T] & \\ \\
\bar{X}_0 =x, \ \bar{Y}_T = - h_x(\bar{X}_T) ,\end{array}\right.
\end{equation}
where $\gamma$ has been defined in (\ref{defdigamma}).
Section 3 is devoted to prove the following result:
\begin{theorem}\label{FBSDE}
Assume {\bf (A)} and  {\bf (B)}, then there exists a unique solution $(\bar{X},\bar{Y},(\bar{Z},\tilde{Z}) \in L^2_\P((0,T)\times \Omega;H)\times L^2_\P((0,T)\times \Omega;D(E))\times L^2_\P((0,T)\times \Omega;L_2(\Xi \times K;H))$ of the forward-backward system \eqref{forwardback_ottimo}. Moreover we have that:
\begin{equation}
\sup_{t \in [0,T[} (T-t)^{1-\alpha}\|\bar{Y_t}\|_{D(E^*)} < + \infty
\end{equation}
\end{theorem}
By the definition of the map $\gamma$ we deduce that $(\gamma([E+ B]^*\bar{Y}),\bar{X},\bar{Y},(\bar{Z},\tilde{Z}))$ is an admissible 4-tuple.
\subsection{Main result}
We can now state the main result of the paper.
\begin{theorem}\label{main}
Assume {\bf(A)} and {\bf(B)}. There exists a unique optimal pair given by the solution of system \eqref{forwardback_ottimo} for the control  problem \eqref{forward} and \eqref{funzionale}.
\end{theorem}
\Dim 
Thanks to theorem \ref{FBSDE} we have an admissible 4-tuple $(\gamma([E+ B]^T\bar{Y}),\bar{X},\bar{Y},(\bar{Z},\tilde{Z}))$ that, by definition of $\gamma$,  verifies condition \eqref{condizione}. So from theorem \ref{condsuff} we deduce the thesis.
\finedim
\section{Proof of theorem \ref{FBSDE}}
System \eqref{forwardback_ottimo} is an infinite dimensional fully coupled forward-backward system. Besides the difficulties typical of the finite dimensional  FBSDEs, see \cite{MaYong}, there are some additional ones due to the presence of unbounded operators. In particular we need to introduce the graph norm of $E$ and prove some crucial estimates with respect to this stronger norm.
Thanks to dissipativity hypotheses \eqref{dissipativita_datofin}, \eqref{dissipativita_costo} and \eqref{dissipativita_gamma}, the more suitable method to get  a solution is the bridge method used in \cite{HuPeng3} whose infinite dimensional extension will be described in next paragraph.
\subsection{The {\em bridge method} applied to an infinite dimensional system }\label{sez-bridge}
This section is devoted to present the {\em bridge method} to solve the Hamiltonian system
(\ref{forwardback_ottimo}), which is a FBSDE in an
infinite dimensional Hilbert space $H$. According to
this method, introduced in \cite{HuPeng3}, in order to solve
a nonlinear fully coupled FBSDE, a linear auxiliary FBSDE is studied
and then making a sort of convex combination between the affine term in this
linear FBSDE and the nonlinear terms in the original FBSDE it is possible to
arrive at the solution of the original FBSDE.

The main difference between the present paper and \cite{HuPeng3} is that
in \cite{HuPeng3} the finite dimensional case is treated and so 
the linear auxiliary FBSDE has a very special structure and is
solvable by hand; in the present paper, since also $Y$ takes its values in $H$,
the auxiliary linear affine FBSDE has a different structure and it takes some efforts to be solved
see section \ref{sez-LQ}.
Namely, let $b_0, h_0 \in  L^2_\P([0,T]\times \Omega; H)$ and $g_0 \in L^2 (\Omega, \F_T;H)$, consider: 
\begin{equation}\label{forwardback_linear_bridge}
\left\{\begin{array}{ll}  d\bar{X}_t = A \bar{X}_t , dt - [E + B][E+ B]^*\bar{Y}_t\,dt+ b_0(t)\,dt + (\lambda -A)D_1\,d\tilde{W}_t + G\, dW_t &\\ \\
 - d{\bar{Y}}_t = A^* {\bar{Y}}_t \, dt   +\bar{X}_t \, dt + h_0(t) \, dt  -\bar{Z}_t \, dW_t  -\tilde{Z}_t \, d\tilde{W}_t , \qquad t \in [0,T] & \\ \\
\bar{X}_0 =x, \, -\bar{Y}_T = \bar{X}_T + g_0 \end{array}\right.
\end{equation}
In the next section we prove the following proposition,
according to which (\ref{forwardback_linear_bridge}) admits a
unique solution. The difficulties in solving this FBSDE
comes at first by the fact that the BSDE contains $Y$ itself, unlike in \cite{HuPeng3},
and by the presence of the unbounded term$[E + B][E+ B]^*\bar{Y}$.
\begin{proposition}\label{prop_forwardbackward_lineare}
Let $b_0, h_0 \in  L^2_\P([0,T]\times \Omega; H)$ and $g_0 \in L^2 (\Omega, \F_T;H)$,
let also $A$, $E$, $B$, $D_1$ and $G$ satisfy assumptions ${\bf(A)}$, then
the linear FBSDE (\ref{forwardback_linear_bridge})
admits a unique mild solution $(\bar X, \bar Y, (\bar Z, \tilde Z))
\in L^2_\P (\Omega;C([0,T];H)) \times L^2_\P (\Omega;C([0,T];H)) \times L^2_\P(\Omega \times [0,T]; L_2(\Xi \times K;H))$ satisfying moreover $$\E\sup_{t\in[0,T]}(T-t)^{2(1-\alpha)}\Vert \bar{Y}_t\Vert_{D((\lambda-A^*)^{1-\alpha})}^2<+\infty.$$
\end{proposition}
The proof of this proposition is given in the next section.

The aim of the present section is to prove the following result on the \emph{bridge method},
in which the solution of the FBSDE (\ref{forwardback_linear_bridge}) is in some sense connected to
the solution of the starting FBSDE (\ref{forwardback_ottimo}).

\noindent Namely, let us define, for $x\in H$, $y\in H\cap D(E)$ and for $\alpha\in[0,1]$,
\begin{align}\label{coeff_fbsde_bridge}
 &b^\alpha(y)=\alpha[E+B]\gamma([E+B]^Ty)+(1-\alpha)[E+B][E+B]^T(-y)\\ \nonumber
&h^\alpha(t,x)=\alpha l^0_x(t,x)+(1-\alpha)(x)\\ \nonumber
&g^\alpha(x)=\alpha h_x(x)-(1-\alpha)(x).
\end{align}
Consider the following FBSDE.
\begin{equation}\label{forwardback_bridge}
\left\{\begin{array}{ll}  d\bar{X}_t = A \bar{X}_t \, dt+b^\alpha(\bar{Y}_t)\,dt+ b_0(t)\,dt + (\lambda -A)D_1\,d\tilde{W}_t + G\, dW_t &\\ \\
 - d{\bar{Y}}_t = A^* {\bar{Y}}_t \, dt   +h^\alpha(\bar{X}_t) \, dt + h_0(t) \, dt  -\bar{Z}_t \, dW_t  -\tilde{Z}_t \, d\tilde{W}_t , \qquad t \in [0,T] & \\ \\
\bar{X}_0 =x, \, -\bar{Y}_T = g^\alpha(\bar X_t)+g_0.\end{array}\right.
\end{equation}
This is, with $\alpha$ varying in $[0,1]$, the systems that links the linear FBSDE
(\ref{forwardback_linear_bridge}) to the original FBSDE (\ref{forwardback_ottimo}).

Notice that, as stated in proposition \ref{prop_forwardbackward_lineare},
the linear FBSDE (\ref{forwardback_linear_bridge}), which
is equal to the FBSDE (\ref{forwardback_bridge}) with $\alpha=0$,
admits an adapted solution satisfying moreover
$\E\sup_{t\in[0,T]}(T-t)^{2(1-\alpha)}\Vert (E+B)\bar{Y}_t\Vert^2<+\infty$. 
In the next lemma we prove that (\ref{forwardback_bridge}) admits a solution.
\begin{lemma}\label{lemma_bridge}
Let $A$, $E$, $B$, $D_1$ and $G$ satisfy assumptions ${\bf(A)}$ and
$\gamma$, $l$ and $h$ satisfy assumptions ${\bf(B)}$.
Assume that for some $\alpha=\alpha_0$ and for any $b_0, h_0 \in  L^2_\P((0,T)\times \Omega; H)$ and any
$g_0 \in L^2 (\Omega, \F_T;H)$ equations (\ref{forwardback_bridge})
admit a mild solution $(\bar X, \bar Y, (\bar Z, \tilde Z))
\in L^2_\P (\Omega;C([0,T];H)) \times L^2_\P (\Omega;C([0,T];H)) 
\times L^2_\P(\Omega \times [0,T]; L_2(\Xi \times K;H))$ satisfying moreover
\[\E\sup_{t\in[0,T]}(T-t)^{2(1-\alpha)}\Vert \bar{Y}_t\Vert_{D((\lambda-A^*)^{1-\alpha})}^2<+\infty.\]
Then there exists $\delta_0 \in (0,1)$
depending only on constants appearing in
${\bf(A)}$ and ${\bf(B)}$
such that for all $\alpha\in[\alpha_0, \alpha_0+\delta]$
and for any $b_0, h_0 \in  L^2_\P([0,T]\times \Omega; H)$ and any
$g_0 \in L^2 (\Omega, \F_T;H)$ FBSDE (\ref{forwardback_bridge})
admits a mild solution $(\bar X, \bar Y, (\bar Z, \tilde Z))
\in L^2_\P (\Omega;C([0,T];H)) \times L^2_\P (\Omega;C([0,T];H)) 
\times L^2_\P(\Omega \times [0,T]; L_2(\Xi \times K;H))$ satisfying moreover
\[\E\sup_{t\in[0,T]}(T-t)^{2(1-\alpha)}\Vert \bar{Y}_t\Vert_{D((\lambda-A^*)^{1-\alpha})}^2<+\infty.\]
\end{lemma}
\Dim
We notice that for $\alpha=\alpha_0+\delta$ coefficients in (\ref{coeff_fbsde_bridge})
can be rewritten as
\begin{align*}
 &b^{\alpha_0+\delta}(y)=b^{\alpha_0}(y)+\delta[E+B]\gamma([E+B]^Ty)+\delta [E+B][E+B]^Ty\\ \nonumber
&h^{\alpha_0+\delta}(t,x)=h^{\alpha_0}(t,x)+\delta l^0_x(t,x)-\delta x\\ \nonumber
&g^{\alpha_0+\delta}(x)=g^{\alpha_0}(x)+\delta h_x(x)+\delta x.
\end{align*}
Notice that by our assumptions it follows that for $\alpha=\alpha_0$
the FBSDE (\ref{forwardback_bridge}) admits a mild solution.

\noindent From this we start proving
that there exists $\delta_0\in (0,1)$ such that for all $\delta \in [0,\delta_0]$,
for all $\alpha\in[\alpha_0,\alpha_0+\delta_0]$ and for all 
$b_0, h_0 \in  L^2_\P([0,T]\times \Omega; H)$ and
 $g_0 \in L^2 (\Omega, \F_T;H)$ the FBSDE (\ref{forwardback_bridge})
admits a unique mild solution $(\bar X, \bar Y, (\bar Z, \tilde Z))
 \in L^2_\P (\Omega;C([0,T];H)) \times L^2_\P (\Omega;C([0,T];H)) 
 \times L^2_\P(\Omega \times [0,T]; L_2(\Xi \times K;H))$ satisfying moreover
 \[ E\sup_{t\in[0,T]}(T-t)^{2(1-\alpha)}\Vert \bar{Y}_t\Vert_{D((\lambda-A^*)^{1-\alpha})}^2<\infty.\]
We set $(\bar X^{0}, \bar Y^{0}, (\bar Z^{0}, \tilde Z^0))=(0,0,(0,0))$.
For $j\geq 0$ we solve iteratively the following FBSDEs
\begin{equation}\label{forwardbackward-iterato}
\left\{\begin{array}{ll}  d\bar{X}^{j+1}_t = A \bar{X}^{j+1}_t \, dt
+\left(\alpha_0[E+B]\gamma([E+B]^* \bar{Y}_t^{j+1})-(1-\alpha_0)[E+B][E+B]^*\bar{Y}_t^{j+1}\right)\,dt
& \\ \\
 \qquad\qquad+\left(\delta[E+B]\gamma([E+B]^* \bar{Y}_t^{k,j})+\delta[E+B][E+B]^*\bar{Y}_t^{j}\right)\,dt & 
\\ \\ \qquad\qquad+b_0(t)\,dt+(\lambda -A)D_1\,d\tilde{W}_t +  G\, dW_t& \\ \\
 - d{\bar{Y}}^{j+1}_t = A_n^* \bar{Y}^{j+1}_t \, dt 
 + h_0(t) \, dt +\alpha_0 l^0_x(t,\bar{X}^{j+1}_t) \, dt+(1-\alpha_0)\bar{X}^{j+1}_t\,dt& \\ \\
 \qquad\qquad\;
+\delta (l^0_x(t,\bar{X}^{j}_t) -\bar{X}^{j}_t)\,dt
 -\bar{Z}^{j+1}_t \, dW_t  -\tilde{Z}^{j+1}_t \, d\tilde{W}_t , \qquad t \in [0,T] & \\ \\
\bar{X}^{j+1}_0 =x, \, -\bar{Y}^{j+1}_T = g^{\alpha_0}(\bar{X}^{j+1}_T)+\delta h_x(\bar{X}^{j}_T)
+\delta \bar{X}^{j}_T+g_0.\end{array}\right.
\end{equation}
Notice that by induction, and by generalizing some statements from the strong to the mild formulation,
see few lines below, such a FBSDE admits a mild solution
$(\bar X^{j+1}, \bar Y^{j+1}, (\bar Z^{j+1}, \tilde Z^{j+1}))
 \in L^2_\P (\Omega;C([0,T];H)) \times L^2_\P (\Omega;C([0,T];H)) 
 \times L^2_\P(\Omega \times [0,T]; L_2(\Xi \times K;H))$ satisfying moreover
 \[\E\sup_{t\in[0,T]}(T-t)^{2(1-\alpha)}\Vert \bar{Y}^{j+1}_t\Vert_{D((\lambda-A^*)^{1-\alpha})}^2<\infty.\]
Indeed, for $j=1$, FBSDE (\ref{forwardbackward-iterato}) is equal to FBSDE (\ref{forwardback_bridge}).
So by hypothesis the solution, with the required regularity, exists.
By induction, assume that for $j$ a solution, with the required regularity exists, and we show that
also for $j+1$ a solution exists. By setting $\tilde{b}_0(t)=\delta[E+B]\gamma([E+B]^* \bar{Y}_t^{j})
+\delta[E+B][E+B]^*\bar{Y}_t^{j}+b_0(t)$; $\tilde{h}_0(t)=\delta (l^0_x(t,\bar{X}^{j}_t) -\bar{X}^{n,j}_t)+h_0(t)$;
$\tilde{g}_0=\delta h_x(\bar{X}^{j}_T)
+\delta \bar{X}^{j}_T+g_0$, FBSDE (\ref{forwardbackward-iterato}) is equal to FBSDE (\ref{forwardback_bridge})
with $\tilde{b}_0$, $\tilde{h}_0$ and $\tilde{g}_0$ in the place of $b_0$, $h_0$ and $g_0$ respectively.
This time $\tilde{b}_0 \notin L^2_\P(\Omega \times [0,T]; H)$, indeed $\tilde{b}_0$
is not well defined as an element of $H$.
Nevertheless, in the mild formulation of $X_t$, $\tilde{b}_0$ appears in integral form
and it is affected by the regularizing properties of the semigroup:
the integral $\int_t^Te^{(s-t)A}\tilde{b}_0(s) ds$
is well defined and bounded in $ L^2_\P(\Omega \times [0,T]; H)\cap L^2_\P(\Omega,C( [0,T]; H))$.

So by our assumptions, $\forall j\geq 0$, there exists a mild solution
$(\bar X^{j+1}, \bar Y^{j+1}, (\bar Z^{j+1}, \tilde Z^{j+1}))
 \in L^2_\P (\Omega;C([0,T];H)) \times L^2_\P (\Omega;C([0,T];H)) 
 \times L^2_\P(\Omega \times [0,T]; L_2(\Xi \times K;H))$ satisfying moreover
 $\E\sup_{t\in[0,T]}(T-t)^{2(1-\alpha)}
\Vert \bar{Y}^{j+1}_t\Vert_{D((\lambda-A)^{1-\alpha})}^2<\infty$.

Next we define, for every $t\in [0,T]$,
\[
\hat{X}^{j+1}_t=\bar{X}^{j+1}_t-\bar{X}^{j}_t;\quad 
\hat{Y}^{j+1}_t=\bar{Y}^{j+1}_t-\bar{Y}^{j}_t;
\]
\[\hat{Z}^{j+1}_t=\bar{Z}^{j+1}_t-\bar{Z}^{j}_t;\quad 
\hat{\tilde Z}^{j+1}_t=\tilde{Z}^{j+1}_t-\tilde{Z}^{j}_t.
\]
We note that $(\hat{X}^{j+1},\hat{Y}^{j+1},(\hat{Z}^{j+1},\hat{\tilde Z}^{j+1}))$
solve
\begin{equation*}
\left\{\begin{array}{ll}  d\hat{X}^{j+1}_t = A \hat{X}^{j+1}_t \, dt
+\alpha_0[E+B](\gamma([E+B]^* \bar{Y}_t^{j+1})-\gamma([E+B]^* \bar{Y}_t^{j}))\,dt& \\ \\
 \qquad\qquad-(1-\alpha_0)[E+B][E+B]^*\hat{Y}_t^{j+1}\,dt
+\delta[E+B][E+B]^*\hat{Y}_t^{j}\,dt
& \\ \\
 \qquad\qquad+\delta[E+B](\gamma([E+B]^* \bar{Y}_t^{j})
-\gamma([E+B]^* \bar{Y}_t^{j-1}))\,dt,& \\ \\
 - d{\hat{Y}}^{j+1}_t = A^* \hat{Y}^{j+1}_t \, dt 
  +\alpha_0 (l^0_x(t,\bar{X}^{j+1}_t)-l^0_x(t,\bar{X}^{j}_t))\,dt
+(1-\alpha_0)\hat{X}^{j+1}_t\,dt& \\ \\
\qquad\qquad +\delta(l^0_x(t,\bar{X}^{j}_t)-l^0_x(t,\bar{X}^{j-1}_t))\,dt-\delta\hat{X}^{j}_t\,dt
 -\hat{Z}^{j+1}_t \, dW_t  -\hat{\tilde{Z}}^{j+1}_t \, d\tilde{W}_t , \qquad t \in [0,T] & \\ \\
\hat{X}^{j+1}_0 =x, 
& \\ \\ -\hat{Y}^{j+1}_T = \alpha_0(h_x(\bar{X}^{j+1}_T)-h_x(\bar{X}^{j}_T))
-(1-\alpha_0)\hat{X}^{j+1}_T
+\delta( h_x(\bar{X}^{j}_T)-h_x(\bar{X}^{j-1}_T))+\delta \hat{X}^{j}_T.\end{array}\right.
\end{equation*}
We notice that by our assumption for every $j$,
$\E\sup_{t\in[0,T]}\vert \hat{X}^{j}_t \vert^2<+\infty$.

\noindent Next we have to apply It\^o formula: in order to do this we have to approximate $X$ and $Y$ with
elements of the domain of $A$. Namely, for $n>\lambda$, we denote as usual $R(n,A):= (n-A)^{-1}$.
We set $(\hat{X}^{n,j+1}, \hat{Y}^{n,j+1}, (\hat{Z}^{n,j+1}, \hat{\tilde{Z}}^{n,j+1}))=
(nR(n,A)\hat{X}^{j+1}, nR(n,A)\hat{Y}^{j+1}, nR(n,A)(\hat{Z}^{j+1}, nR(n,A)\hat{\tilde{Z}}^{j+1}))$.
We also denote $E_n+B_n:=nR(n,A)(E+B)$ and we note that $(\hat{X}^{n,j+1}, \hat{Y}^{n,j+1}, (\hat{Z}^{n,j+1}, \hat{\tilde{Z}}^{n,j+1}))$ solve the following
\begin{equation*}
\left\{\begin{array}{ll}  d\hat{X}^{n,j+1}_t = A \hat{X}^{n,j+1}_t \, dt
+\alpha_0[E_n+B_n](\gamma([E+B]^T \bar{Y}_t^{j+1})-\gamma([E+B]^T \bar{Y}_t^{j}))\,dt& \\ \\
 \qquad\qquad-(1-\alpha_0)[E_n+B_n][E+B]^*\hat{Y}_t^{j+1}\,dt
+\delta[E_n+B_n][E+B]^*\hat{Y}_t^{j}\,dt
& \\ \\
 \qquad\qquad+\delta[E_n+B_n](\gamma([E+B]^* \bar{Y}_t^{j})
-\gamma([E+B]^* \bar{Y}_t^{j-1}))\,dt,& \\ \\
 - d{\hat{Y}}^{n,j+1}_t = A^* \hat{Y}^{n,j+1}_t \, dt 
  +\alpha_0 nR(n,A)(l^0_x(t,\bar{X}^{j+1}_t)-l^0_x(t,\bar{X}^{j}_t))\, dt+(1-\alpha_0)\hat{X}^{n,j+1}_t\, dt& \\ \\
\qquad\qquad +\delta nR(n,A)(l^0_x(t,\bar{X}^{j}_t)-l^0_x(t,\bar{X}^{j-1}_t))\, dt-\delta\hat{X}^{n,j}_t\, dt
 -\hat{Z}^{n,j+1}_t \, dW_t  -\hat{\tilde{Z}}^{n,j+1}_t \, d\tilde{W}_t , \qquad t \in [0,T] & \\ \\
\hat{X}^{n,j+1}_0 =x, 
& \\ \\ -\hat{Y}^{n,j+1}_T = \alpha_0 nR(n,A)(h_x(\bar{X}^{j+1}_T)-h_x(\bar{X}^{j}_T))
-(1-\alpha_0)\hat{X}^{n,j+1}_T
& \\ \\ \qquad\qquad +\delta nR(n,A)( h_x(\bar{X}^{j}_T)-h_x(\bar{X}^{j-1}_T))+\delta \hat{X}^{n,j}_T.\end{array}\right.
\end{equation*}
By applying It\^o formula to $\langle \hat{X}^{n,j+1}_t,\hat{Y}^{n,j+1}_t \rangle$, and then integrating over $[0,T]$ and taking
expectation we get
\begin{align}\label{bridge-ito-n}
 &-\E\< \hat{X}^{n,j+1}_T,\alpha_0 nR(n,A)(h_x(\bar{X}^{j+1}_T)-h_x(\bar{X}^{j}_T))
-(1-\alpha_0)\hat{X}^{n,j+1}_T\\
&-\delta\< nR(n,A)(h_x(\bar{X}^{j}_T)-h_x(\bar{X}^{j-1}_T))+ \hat{X}^{n,j}_T,\hat{X}^{n,j+1}_T\> \\ \nonumber
&=\E\int_0^T\left[ \< \alpha_0[E_n+B_n](\gamma([E+B]^* \bar{Y}_t^{j+1})-
\gamma([E+B]^* \bar{Y}_t^{j})),\hat{Y}_t^{n,j+1}\> \right.\\ \nonumber
&\left.+(1-\alpha_0)\< [E_n+B_n][E+B]\hat{Y}_t^{j+1},\hat{Y}_t^{n,j+1}\vert^2 \right] dt \\ \nonumber
&-\E\int_0^T\left[\alpha_0\< nR(n,A)(l^0_x(t,\bar{X}^{j+1}_t)-l^0_x(t,\bar{X}^{j}_t)),\hat{X}^{n,j+1}_t \>
+(1-\alpha_0)\vert \hat{X}^{n,j+1}_t \vert^2\right]\,dt\\ \nonumber
&+\delta\E \int_0^T \< (\gamma([E+B]^* \bar{Y}_t^{j})-\gamma([E+B]^T \bar{Y}_t^{j-1}))+ 
[E+B]^*\hat{Y}_t^{j},[E_n+B_n]^*\hat{Y}_t^{n,j+1}\> \, dt \\  \nonumber
&-\delta\E\int_0^T \< nR(n,A)(l^0_x(t,\bar{X}^{j}_t)-l^0_x(t,\bar{X}^{j-1}_t))+\hat{X}^{n,j}_t ,\hat{X}^{n,j+1}_t \> \,dt. \nonumber
\end{align}
Next we want to let $n\rightarrow +\infty$ in the (\ref{bridge-ito-n}):
in order to do this we need to recover at first 
the $L^2_\P(\Omega \times [0,T]; H)$-convergence of $\hat{X}^{n,j}$
and of $\hat{Y}^{n,j}$ to $\hat{X}^{j}$ and of $\hat{Y}^{j}$ respectively.
For what concerns $\hat{X}^{n,j}$, notice that in mild form
\begin{align*}
 \hat{X}^{n,j+1}_t- \hat{X}^{j+1}_t&=\int_0^te^{(t-s) A} 
\alpha_0(nR(n,A)-I)[E+B](\gamma([E+B]^* \bar{Y}_t^{j+1})-\gamma([E+B]^* \bar{Y}_t^{j}))\,dt \\
&-\int_0^te^{(t-s) A}(1-\alpha_0)(nR(n,A)-I)[E+B][E+B]^*\hat{Y}_t^{j+1}\,dt\\
&+\delta\int_0^te^{(t-s) A}(nR(n,A)-I)[E+B][E+B]^*\hat{Y}_t^{j}\,dt \\
&+\delta\int_0^t e^{(t-s) A}(nR(n,A)-I)[E+B](\gamma([E+B]^* \bar{Y}_t^{j})
-\gamma([E+B]^* \bar{Y}_t^{j-1}))\,dt.
\end{align*}
Since $e^{(t-s) A}$ and $(nR(n,A)-I)$ commute, we get
\begin{align*}
 \hat{X}^{n,j+1}_t- \hat{X}^{j+1}_t&=\int_0^t (nR(n,A)-I) e^{(t-s) A} 
\alpha_0[E+B](\gamma([E+B]^* \bar{Y}_t^{j+1})-\gamma([E+B]^* \bar{Y}_t^{j}))\,dt \\
&-\int_0^t (nR(n,A)-I) e^{(t-s) A}(1-\alpha_0)[E+B][E+B]^*\hat{Y}_t^{j+1}\,dt\\
&+\delta\int_0^t (nR(n,A)-I) e^{(t-s) A}[E+B][E+B]^*\hat{Y}_t^{j}\,dt \\
&+\delta\int_0^t(nR(n,A)-I)e^{(t-s) A}[E+B](\gamma([E+B]^* \bar{Y}_t^{j})
-\gamma([E+B]^* \bar{Y}_t^{j-1}))\,dt.
\end{align*}
Let us consider the first integral, for the others the same conclusion follows in a similar way.
For a.a. $s\in[0,t]$ and for all $0< t\leq T$, and for a.a. $\omega \in \Omega$,
$\alpha_0 e^{(t-s) A} 
[E+B](\gamma([E+B]^* \bar{Y}_t^{j+1})-\gamma([E+B]^* \bar{Y}_t^{j}))$ is
an element in $H$ and 
$$\vert \alpha_0 e^{(t-s) A} 
[E+B](\gamma([E+B]^* \bar{Y}_t^{j+1})-\gamma([E+B]^* \bar{Y}_t^{j})) \vert\leq C\Delta \alpha_0 (t-s)^{-(1-\alpha)}\vert
[E+B]^* \hat{Y}_t^{j+1}\vert,$$
so that, by dominated convergence, as $n\rightarrow \infty$:
\[
\E\sup_{t\in[0,T]}\vert \hat{X}_t^{n,j+1}- \hat{X}_t^{j+1}\vert^2 \rightarrow 0.
\]

\noindent In a similar way we can get that $\hat{Y}_t^{n,j+1}\rightarrow \hat{Y}_t^{j+1}$
in $L^2_\P(\Omega\times[0,T];H)$ and moreover  as $n\rightarrow \infty $:
 \[\E\sup_{t\in[0,T]}(T-t)^{2(1-\alpha)} \Vert (\hat{Y}_t^{n,j+1}- \hat{Y}_t^{j+1})\Vert
_{D((\lambda-A^*)^{1-\alpha})}^2\rightarrow 0. \]

\noindent In order to let $n\rightarrow \infty$ in (\ref{bridge-ito-n}) we have also to show that
$\E\sup_{t\in[0,T]}(T-t)^{2(1-\alpha)}
\vert [E_n+B_n]^*\hat{Y}_t^{n,j+1}- [E+B]^*\hat{Y}_t^{j+1})\vert^2\rightarrow 0 $ as $n\rightarrow \infty$.
Notice that,
\begin{align*}
&[E_n+B_n]^*\hat{Y}^{n,j+1}_t-[E+B]^*\hat{Y}^{j+1}_t\\
& =-\Et e^{(T-t) A^*}  \left(
\alpha_0([E_n+B_n]^*nR(n,A)-[E+B]^*)(h_x(\bar{X}^{j+1}_T)-h_x(\bar{X}^{j}_T))\right.\\
&\left.-(1-\alpha_0)(nR(n,A)-I)\hat{X}^{j+1}_T
+\delta (nR(n,A)-I)( h_x(\bar{X}^{j}_T)-h_x(\bar{X}^{j-1}_T))+\delta (nR(n,A)-I)\hat{X}^{j}_T\right)\\
& + \Et\int_t^T ([E_n+B_n]^*nR(n,A)-[E+B]^*)[\alpha_0(l^0_x(s,\bar{X}^{j+1}_s)-l^0_x(s,\bar{X}^{j}_s))
+(1-\alpha_0)\hat{X}^{j+1}_s]\,ds \\
& +\Et\int_t^T  \delta([E_n+B_n]^*nR(n,A)-[E+B]^*)(l^0_x(s,\bar{X}^{j}_s)-l^0_x(s,\bar{X}^{j-1}_s)-\delta\hat{X}^{j}_s)\, ds. \\
\end{align*}
So we can let $n\rightarrow\infty$ in (\ref{bridge-ito-n}) and we get
\begin{align*}
 &-\E\< \hat{X}^{j+1}_T,\alpha_0 (h_x(\bar{X}^{j+1}_T)-h_x(\bar{X}^{j}_T))
-(1-\alpha_0)\hat{X}^{j+1}_T\rangle\\
&-\delta\< (h_x(\bar{X}^{j}_T)-h_x(\bar{X}^{j-1}_T))+ \hat{X}^{j}_T,\hat{X}^{j+1}_T\> \\ \nonumber
&=\E\int_0^T\left[ \< \alpha_0[E+B](\gamma([E+B]^* \bar{Y}_t^{j+1})-
\gamma([E+B]^* \bar{Y}_t^{j})),\hat{Y}_t^{j+1}\> \right.\\ \nonumber
&\left.+(1-\alpha_0)\< [E+B][E+B]\hat{Y}_t^{j+1},\hat{Y}_t^{j+1}\vert^2 \right] dt \\ \nonumber
&-\E\int_0^T\left[\alpha_0\< (l^0_x(t,\bar{X}^{j+1}_t)-l^0_x(t,\bar{X}^{j}_t)),\hat{X}^{j+1}_t \>
+(1-\alpha_0)\vert \hat{X}^{j+1}_t \vert^2\right]\,dt\\ \nonumber
&+\delta\E \int_0^T \< (\gamma([E+B]^* \bar{Y}_t^{j})-\gamma([E+B]^T \bar{Y}_t^{j-1}))+ 
[E+B]^*\hat{Y}_t^{j},[E+B]^*\hat{Y}_t^{j+1}\> \, dt \\  \nonumber
&-\delta\E\int_0^T \< (l^0_x(t,\bar{X}^{n,j}_t)-l^0_x(t,\bar{X}^{j-1}_t))+\hat{X}^{j}_t ,\hat{X}^{j+1}_t \> \,dt. \nonumber
\end{align*}
So, by assumptions ${\bf(B)}$ we get
\begin{align*}
 & \min {\{c_1,1\} }\E \vert\hat{X}^{j+1}_T\vert^2\\
& \leq \delta\Delta \E\vert\hat{X}^{j+1}_T\vert 
\,\vert\hat{X}^{j}_T\vert
-\alpha_0 c_1\E\int_0^T\vert[E+B]^* \hat{Y}_t^{j+1}\vert^2\, dt-(1-\alpha_0)\E\int_0^T\vert[E+B]^*\hat{Y}_t^{j+1}\vert^2\,  dt \\
&+\delta(\Delta+1)\E \int_0^T \vert[E+B]^* \hat{Y}_t^{j}\vert \vert\hat{Y}_t^{j+1}\vert\, dt -(\alpha_0 c_1+1-\alpha_0)\E\int_0^T\vert \hat{X}^{j+1}_t\vert^2\,dt\\
&+\delta\Delta\E\int_0^T \vert\hat{X}^{j}_t\vert \,\vert\hat{X}^{j+1}_t \vert \,dt.
\end{align*}
By applying Young inequalities several times we finally get
\begin{align*}
 & \E\vert\hat{X}^{j+1}_T\vert^2+[\E\int_0^T\vert[E+B]^*\hat{Y}_t^{j+1}\vert^2  dt +
\E\int_0^T\vert\hat{X}_t^{j+1}\vert^2  dt ])\\
& \leq c'( \delta,\Delta,c_1) \E\vert\hat{X}^{j}_T\vert^2 
+ c'(\delta,c_1)\E\int_0^T\vert[E+B]^*\hat{Y}_t^{j}\vert^2  dt 
+c'(\delta,\Delta,c_1)\E\int_0^T\vert[E+B]^*\hat{Y}_t^{j}\vert^2  dt,
\end{align*}
where $c'( \delta,\Delta,c_1)$ and $c'(\delta,c_1)$
are constants depending respectively only on $\delta\,,\Delta\,,c_1$ and $\delta\,,c_1$
respectively.
Now notice that
\begin{align*}
 &\hat{X}^{j}_T=\alpha_0 \int_0^Te^{(T-t)A} [E+B](\gamma([E+B]^* \bar{Y}_t^{j})-\gamma([E+B]^* \bar{Y}_t^{j-1}))\,dt\\
&+(1-\alpha_0) \int_0^Te^{(T-t)A} [E+B][E+B]^* \hat{Y}_t^{j}\,dt
+\delta\int_0^Te^{(T-t)A} [E+B][E+B]^* \hat{Y}_t^{j-1}\,dt\\
&+\delta \int_0^Te^{(T-t)A} [E+B](\gamma([E+B]^* \bar{Y}_t^{j-1})
-\gamma([E+B]^* \bar{Y}_t^{j-2}))\,dt.
\end{align*}
So
\begin{align*}
 &\E\vert\hat{X}^{j}_T\vert^2\leq \alpha_0^2 \E\vert\int_0^T\Vert e^{(T-t)A} [E+B]\Vert 
\Delta\vert[E+B]^* \hat{Y}_t^{j}\vert\,dt\vert^2\\
&+(1-\alpha_0)^2 \E\vert\int_0^T c(T-t)^{-(1-\alpha)} \vert[E+B]^* \hat{Y}_t^{j}\vert\,dt\vert^2
+\delta^2\E\vert\int_0^T \Vert e^{(T-t)A} [E+B]\Vert\vert[E+B]^* \hat{Y}_t^{j-1})\vert\,dt\vert^2\\
&+\delta^2 \E\vert\int_0^T\Vert e^{(T-t)A} [E+B]\Vert\delta\vert[E+B]^* \hat{Y}_t^{j-1}\vert\,dt\vert^2\\
&\leq c\Delta^2\alpha_0^2 \E\vert\int_0^T(T-t)^{-(1-\alpha)} 
\vert[E+B]^* \hat{Y}_t^{j}\vert\,dt\vert^2+c(1-\alpha_0)^2 \E\int_0^T (T-t)^{-2(1-\alpha)} \int_0^T \vert[E+B]^* \hat{Y}_t^{j}\vert^2\,dt
\\
&+c\delta^2\E\vert\int_0^T (T-t)^{-(1-\alpha)}\vert[E+B]^* \hat{Y}_t^{j-1})\vert\,dt\vert^2
+\delta^2 \E\vert\int_0^T (T-t)^{-(1-\alpha)} [E+B]\Vert\delta\vert[E^k+B]^* \hat{Y}_t^{j-1}\vert\,dt\vert^2\\
&\leq c\Delta^2\alpha_0^2 \E\sup_{t\in[0,T]}(T-t)^{2(1-\alpha)}\vert[E+B]^* \hat{Y}_t^{j}\vert^2
\left(\int_0^T(T-t)^{-2(1-\alpha)}\,dt\right)^2\\
&+c(1-\alpha_0)^2 \E\sup_{t\in[0,T]}(T-t)^{2(1-\alpha)}\vert[E+B]^* \hat{Y}_t^{j}\vert^2\left(\int_0^T (T-t)^{-2(1-\alpha)} \,dt\right)^2\\
&+c\delta^2\E\sup_{t\in[0,T]}(T-t)^{2(1-\alpha)}\vert[E+B]^* \hat{Y}_t^{j}\vert^2\left(\int_0^T (T-t)^{-2(1-\alpha)}\,dt\right)^2\\
&+c\delta^2\E\sup_{t\in[0,T]}(T-t)^{2(1-\alpha)}\vert[E+B]^* \hat{Y}_t^{j}\vert^2
\left(\int_0^T (T-t)^{-2(1-\alpha)}\,dt\right)^2.
\end{align*}
Now, arguing as in \cite{HuPeng3}, proof of lemma 3.2, we get that
there exists $\delta_0 \in (0,1)$ depending only on $c_1, \Delta, T$, such that
for every $\delta\in (0,\delta_0]$, we get
\begin{align*}
&\E\int_0^T\vert[E+B]^*\hat{Y}_t^{j+1}\vert^2  dt +
\E\int_0^T\vert\hat{X}_t^{j+1}\vert^2  dt ]\\
 &\leq \dfrac{1}{4}\left[ \E\int_0^T\vert\hat{X}_t^{j}\vert^2  dt
+ \E\int_0^T\vert[E+B]^*\hat{Y}_t^{j}\vert^2  dt \right]\\
&+\dfrac{1}{8}\left[ \E\int_0^T\vert\hat{X}_t^{j-1}\vert^2  dt
+ \E\int_0^T\vert[E+B]^*\hat{Y}_t^{j-1}\vert^2  dt \right].
\end{align*}
From this we deduce that
$(\bar{X}_t^{j},\bar{Y}_t^{j})_{j\geq 1}$ is a Cauchy sequence in $L^2_{\mathcal P}(\Omega\times[0,T],H)
\times L^2(\Omega\times[0,T],H)$ and we denote by $(\bar{X}_t,\bar{Y}_t)$ its limit.

\noindent In order to prove that $(\bar{X}_t^{j},\bar{Y}_t^{j})_{j\geq 1}$
converge to $(\bar{X}_t,\bar{Y}_t)$ also in $L^2_{\mathcal P}(\Omega,C([0,T],H))
\times L^2_{\mathcal P}(\Omega,C([0,T],H))$ we go to the mild formulation of the equations
solved by $\bar{X}_t^{j}$ and $\bar{Y}_t^{j}$. We start by $\bar{X}_t^{j}$:
\begin{align*}
 &\bar{X}^{j}_t=e^{t A}x+\alpha_0 \int_0^te^{(t-s)A} [E+B](\gamma([E+B]^* \bar{Y}_s^{j})\,ds\\
&-(1-\alpha_0) \int_0^te^{(t-s)A} [E+B][E+B]^* \bar{Y}_s^{j}\,ds
+\delta\int_0^te^{(t-s)A} [E+B][E+B]^* \bar{Y}_s^{j-1}\,ds\\
&+\delta \int_0^Te^{(t-s)A} [E+B](\gamma([E+B]^* \bar{Y}_s^{j-1})\,ds+\int_0^te^{(t-s)A}b_0(s)\,ds\\
&+\int_0^t e^{(t-s)A}G\,dW_s+\int_0^t e^{(t-s)A}D_1\,d\tilde{W}_s.
\end{align*}
So
\begin{align*}
 &\E\sup_{t\in[0,T]}\vert\bar{X}^{j}_t\vert^2\\
&\leq \alpha_0^2 \E\sup_{t\in[0,T]}\vert\int_0^t\Vert e^{(t-s)A} [E+B]\Vert 
\Delta\vert[E+B]^* \bar{Y}_s^{j}\vert\,ds\vert^2\\
&+(1-\alpha_0)^2\E\sup_{t\in[0,T]}(T-t)^{2(1-\alpha)}\vert[E+B]^* \bar{Y}_t^{j}\vert^2 
\sup_{t\in[0,T]}\left(\int_0^t c(t-s)^{-(1-\alpha)} (T-s)^{-(1-\alpha)}\,ds\right)^2 \\
&+\delta^2\E\vert\sup_{t\in[0,T]}\int_0^t \Vert e^{(t-s)A} [E+B]\Vert\vert[E+B]^* 
\bar{Y}_s^{j-1})\vert\,ds\vert^2\\
&+\delta^2 \sup_{t\in[0,T]}\E\vert\int_0^t\Vert e^{(t-s)A} [E+B]\Vert\delta\vert[E+B]^* \bar{Y}_s^{j-1}\vert\,ds\vert^2\\
&\leq c\Delta^2\alpha_0^2 \E\sup_{t\in[0,T]}(T-t)^{2(1-\alpha)}\vert[E+B]^* \bar{Y}_t^{j}\vert^2 
\sup_{t\in[0,T]}\left(\int_0^t c(t-s)^{-(1-\alpha)} (T-s)^{-(1-\alpha)}\,ds\right)^2 \\
&+c(1-\alpha_0)^2\E\sup_{t\in[0,T]}(T-t)^{2(1-\alpha)}\vert[E+B]^* \bar{Y}_t^{j}\vert^2 
\sup_{t\in[0,T]}\left(\int_0^t c(t-s)^{-(1-\alpha)} (T-s)^{-(1-\alpha)}\,ds\right)^2 
\\
&+c\delta^2\E\sup_{t\in[0,T]}(T-t)^{2(1-\alpha)}\vert[E+B]^* \bar{Y}_t^{j-1}\vert^2 
\sup_{t\in[0,T]}\left(\int_0^t c(t-s)^{-(1-\alpha)} (T-s)^{-(1-\alpha)}\,ds\right)^2 \\
&+\delta^2 \E\sup_{t\in[0,T]}(T-t)^{2(1-\alpha)}\vert[E+B]^* \bar{Y}_t^{j-1}\vert^2 
\sup_{t\in[0,T]}\left(\int_0^t c(t-s)^{-(1-\alpha)} (T-s)^{-(1-\alpha)}\,ds\right)^2 \\
\end{align*}
By the previous choice of $\delta$ we get that $(\bar{X}^{j})_{j\geq 1}$ is a Cauchy sequence
in $L^2_{\mathcal P}(\Omega,C([0,T],H))$ so that $\bar{X}^{j}\rightarrow \bar{X}$ in
$L^2_{\mathcal P}(\Omega,C([0,T],H))$. For what concerns the convergence of $\bar{Y}^{j}$
in $L^2_{\mathcal P}(\Omega,C([0,T],H))$ we have first to recover the convergence of
$(\bar{Z}^{j}, \tilde{Z}^{j})$ in $L^2_\P(\Omega \times [0,T]; L_2(\Xi \times K;H))$.
In its mild formulation, $\bar{Y}^{j+1}$ solves the following BSDE
\begin{align}\label{bsde-mild}
 \bar{Y}^{j+1}_t&=-e^{(T-t) A^*}\left[ \alpha_0h_x(\bar{X}^{j+1}_T)
-(1-\alpha_0)\bar{X}^{j+1}_T
+\delta h_x(\bar{X}^{j}_T)+\delta \bar{X}^{j}_T +g_0\right]\\ \nonumber
&+\alpha_0 \int_t^T e^{(s-t)A^*}l_x^0(s,\bar{X}_s^{j+1})\,ds
+(1-\alpha_0)\int_t^T e^{(s-t)A^*}\bar{X}_s^{j+1}\,ds \\ \nonumber
&+\delta\int_t^T e^{(s-t)A^*} l_x^0(s,\bar{X}_s^{j})\,ds
-\delta \int_t^T e^{(s-t)A^*}\bar{X}_s^{j}\,ds\\ \nonumber
&-\int_t^T e^{(s-t)A^*}\bar{Z}_s^{j+1}\,dW_s
-\int_t^T e^{(s-t)A^*}\tilde{Z}_s^{j+1}\,d\tilde{W}_s +\int_t^T e^{(s-t)A^*}h_0(s)\,ds. \nonumber
\end{align}
Let us denote by 
\begin{align*}
f^{j+1}_s&:=\alpha_0 l_x^0(s,\bar{X}_s^{j+1})
+(1-\alpha_0) e^{(s-t)A^*}\bar{X}_s^{j+1} \\ 
&+\delta  l_x^0(s,\bar{X}_s^{j})\,ds
-\delta  \bar{X}_s+h_0(s).\\
\end{align*}
Arguing as in \cite{HuPeng}, by the extended martingale representation theorem, (see also
\cite{HuPeng1} and \cite{Yor})), for every $s\in[0,T]$ there exists $(K^{j}(s,\cdot),\tilde{K}^{j}(s,\cdot)
\in L^2_\P(\Omega\times[0,T],L_2(\Xi,H))\times L_\P^2(\Omega\times[0,T],L_2(K,H))$
such that $\forall\, 0\leq t\leq s\leq T$
\[
 \Et f^{j}_s=\E f^{j}_s+\int_0^t K^{j}(s,\theta )\,dW_\theta+
\int_0^t K^{j}(s,\theta )\,d\tilde{W}_\theta.
\]
Note that $\forall \,\theta\geq s$, $K^{j}(s,\theta)=0$ a.e.
and $\tilde{K}^{j}(s,\theta)=0$ a.e..; and
\begin{equation}
\label{stimaK}
\E \int_0^T\int_0^s\left[\vert K^{j}(s,\theta)\vert^2
+\vert \tilde{K}^{j}(s,\theta)\vert^2 \right]\, d\theta \,ds
\leq 4\E \int_0^T \vert f^{j}_s \vert^2 \,ds.
\end{equation}
Moreover, there exists $(L^{j},\tilde{L}^{j})\in 
L^2_\P(\Omega\times[0,T],L_2(\Xi,H))\times L^2_\P(\Omega\times[0,T],L_2(K,H))$ such that 
\[
 \Et \bar{Y}^{j+1}_T=\E\bar{Y}^{j+1}_T+\int_0^t L^{j}(\theta )\,dW\theta+
\int_0^t \tilde{L}^{j}(\theta )\,d\tilde{W}_\theta.
\]
So we get
\begin{align*}
 \bar{Y}^{j+1}_t&=e^{(T-t) A^T}\left[ \alpha_0h_x(\bar{X}^{j+1}_T)
-(1-\alpha_0)\bar{X}^{j+1}_T
+\delta h_x(\bar{X}^{j}_T)+\delta \bar{X}^{j}_T \right]\\
&+\alpha_0 \int_t^T e^{(s-t)A^T} f_s^{j+1}\,ds
-\int_t^T e^{(s-t)A^*}L_s^{j+1}\,dW_s
-\int_t^T e^{(s-t)A^*}\tilde{L}_s^{j+1}d\tilde{W}_s \\
&-\int_t^T e^{(s-t)A^*}\int_s^Te^{(\alpha-s)A^*} K^{j}(\alpha,s)\,d\alpha\,dW_s\\
&-\int_t^T e^{(s-t)A^*}\int_s^Te^{(\alpha-s)A^*} \tilde{K}^{j}(\alpha,s)\,d\alpha\,d\tilde{W}_s.
\end{align*}
By comparing with (\ref{bsde-mild}) we deduce that, for almost all $s\in[0,T]$,
\begin{align*}
&\bar{Z}_s^{j}=\int_s^Te^{(\alpha-s)A^T} K_s^{j}(\alpha,s)\,d\alpha,\\
&\tilde{Z}_s^{j}=\int_s^Te^{(\alpha-s)A^T} \tilde{K}_s^{j}(\alpha,s)\,d\alpha.
\end{align*}
By the definition of $(K^{j},\tilde{K}^{j})$, by estimates (\ref{stimaK}), and by previous estimates on
the $L^2$-norm of $\hat{X}^{j}$ and of $\hat{Y}^{j}$ it is possible to prove that
$(\bar{Z}_s^{j}, \tilde{Z}_s^{j})$ is a Cauchy sequence in
$L^2_\P(\Omega \times [0,T]; L_2(\Xi \times K;H))$, and we denote by
$(\bar{Z}_s, \tilde{Z}_s)$ its limit.

\noindent We are ready to prove that $\bar{Y}^{j}\rightarrow
\bar{Y}$ in $L^2_\P(\Omega,C([0,T],H))$. We can rewrite (\ref{bsde-mild}) as
\begin{align*}
 &\Et\bar{Y}^{j+1}_t=\bar{Y}^{j+1}_t\\
&=\Et e^{(T-t) A^*}\left[ -\alpha_0h_x(\bar{X}^{j+1}_T)
-(1-\alpha_0)\bar{X}^{j+1}_T
+\delta h_x(\bar{X}^{j}_T)+\delta \bar{X}^{j}_T \right]\\
&+\alpha_0 \Et\int_t^T e^{(s-t)A^*}l_x^0(s,\bar{X}_s^{j+1})\,ds
+(1-\alpha_0)\Et\int_t^T e^{(s-t)A^*}\bar{X}_s^{j+1}\,ds \\ 
&+\delta\Et\int_t^T e^{(s-t)A^*} l_x^0(s,\bar{X}_s^{j})\,ds
-\delta \Et\int_t^T e^{(s-t)A^*}\bar{X}_s^{j}\,ds\\
&+\E\int_t^T e^{(s-t)A^*}h_0(s)\,ds,
\end{align*}
so that
\begin{align*}
 \E\sup_{t\in[0,T]}\vert\hat{Y}^{j+1}_t\vert^2&\leq c(T,A,\Delta)\left[ \vert\hat{X}^{j+1}_T\vert^2
+\delta^2(1+\Delta^2) \vert\hat{X}^{j}_T\vert^2 \right]\\
&+\alpha_0^2 c(T,A,\Delta)^2 \E\int_0^T \vert\hat{X}_s^{j+1}\vert^2\,ds
+\delta^2c(T,A,\Delta)\E\int_0^T \vert\hat{X}_s^{j}\vert^2\,ds.
\end{align*}
From this, again by using the previous estimates on
the $L^2$-norm of $\hat{X}^{j}$ and of $\hat{Y}^{j}$  it is possible to prove that
$\hat{Y}^{j}$ is a Cauchy sequence in $L^2_\P(\Omega,C([0,T],H))$ and the claim follows.

\noindent Finally we have to prove that $\E\sup_{t\in[0,T]}(T-t)^{2(1-\alpha)}
\Vert \bar{Y}_t \Vert_{D((\lambda-A)^{1-\alpha})}^2 $ is bounded. Let $\alpha\in[\alpha_0,
\alpha_0+\delta]$.
In its mild formulation, $\bar{Y}$ solves the following BSDE
\begin{align*}
 \bar{Y}_t&=e^{(T-t) A^*} \left[-\alpha h_x(\bar{X}_T)
+(1-\alpha)\bar{X}_T+g_0(t)\right]\\
 &+\alpha \int_t^T e^{(s-t)A^*}\left(\alpha l_x^0(s,\bar{X}_s)
 +(1-\alpha)\bar{X}_s\right)\,ds +\int_t^T e^{(s-t)A^*}h_0(s)\,ds\\ 
 &-\int_t^T e^{(s-t)A^*}\bar{Z}_s\,dW_s
 -\int_t^T e^{(s-t)A^*}\tilde{Z}_s\,d\tilde{W}_s . \nonumber
\end{align*}
Notice also that
\begin{align*}
  \bar{Y}_t&=\Et \bar{Y}_t=\Et e^{(T-t) A^*} \left[-\alpha h_x(\bar{X}_T)
  +(1-\alpha)\bar{X}_T+g_0(T)\right]\\
  &+ \Et \int_t^Te^{(s-t)A^*}\left(\alpha l_x^0(s,\bar{X}_s)
  +(1-\alpha)\bar{X}_s\right)\,ds +\Et\int_t^T e^{(s-t)A^*}h_0(s)\,ds
 \end{align*}
By the regularizing properties of the semigroup $(e^{tA})_{t\geq 0}$, by the assumptions on $E$
and by the previous mild equality satisfied by $\bar{Y}_t$ we get that for every
$t\in[0,T]$, $\bar{Y}_t\in \cald(E)$ and 
\begin{align*}
&\E\sup_{t\in[0,T]} \vert(T-t)^{2(1-\alpha)}(E+B)^* \bar{Y}_t\vert^2\\
&\leq c\E \sup_{t\in[0,T]} (T-t)^{2(1-\alpha)}\vert(E+B)^*e^{(T-t) A^*} \left[-\alpha h_x(\bar{X}_T)
  +(1-\alpha)\bar{X}_T+g_0(T)\right]\vert^2\\
  &+ c\E \sup_{t\in[0,T]} (T-t)^{2(1-\alpha)}\vert(E+B)^*
\int_t^Te^{(s-t)A^*}\left(\alpha l_x^0(s,\bar{X}_s)
  +(1-\alpha)\bar{X}_s\right)\,ds\vert^2 \\
&+c\E \sup_{t\in[0,T]} (T-t)^{2(1-\alpha)}\vert(E+B)^*\int_t^T e^{(s-t)A^*}h_0(s)\,ds\vert^2=I+II+III.
 \end{align*}
Recall that $E=(\lambda-A)D$, and $D$ takes its values in $D(\lambda-A)^\alpha$, so that also by the
analyticity of $A$, we get, for every $t>0$ and every $f\in H$
$$
\vert e^{tA} f\vert\leq c t^{-(1-\alpha)}\vert f\vert.
$$
So
\begin{align*}
&I\leq c\E \sup_{t\in[0,T]} (1+\vert\bar{X}_T\vert)^2<+\infty; \\
&II \leq  c\E \sup_{t\in[0,T]} (T-t)^{2(1-\alpha)}
\int_t^T (s-t)^{-(1-\alpha)}\left(1+\vert \bar{X}_s\vert\right)\,ds\vert^2\\
&\leq \E \sup_{t\in[0,T]} (T-t)^{2(1-\alpha)}
\int_t^T (s-t)^{-2(1-\alpha)}\,ds \int_0^T\left(1+\vert \bar{X}_s\vert^2\right)\,ds<+\infty;\\
&III\leq c\E \sup_{t\in[0,T]} (T-t)^{2(1-\alpha)}
\vert\int_t^T\vert (s-t)^{-(1-\alpha)}h_0(s)\,ds\vert^2<+\infty.
 \end{align*}
In order to conclude the proof, note also that $(\bar{X}, \bar{Y},(\bar{Z}, \tilde{Z}))$
is a solution to the FBSDE (\ref{forwardback_bridge}).
\finedim
\begin{remark}\label{remarkBnullo}
 We notice that the presence of a diffuse control is not required in our methods. Indeed, if $B=0$
as an auxiliary linear FBSDE we can consider
\begin{equation*}
\left\{\begin{array}{ll}  d\bar{X}_t = A \bar{X}_t , dt - [E + I][E+I]^*\bar{Y}_t\,dt+ b_0(t)\,dt + (\lambda -A)D_1\,d\tilde{W}_t + G\, dW_t &\\ \\
 - d{\bar{Y}}_t = A^* {\bar{Y}}_t \, dt   +\bar{X}_t \, dt + h_0(t) \, dt  -\bar{Z}_t \, dW_t  -\tilde{Z}_t \, d\tilde{W}_t , \qquad t \in [0,T] & \\ \\
\bar{X}_0 =x, \, -\bar{Y}_T = \bar{X}_T + g_0 \end{array}\right.
\end{equation*}
and we can apply the bridge method linking this FBSDE to the FBSDE
\begin{equation*}
\left\{\begin{array}{ll}  d\bar{X}_t = A \bar{X}_t , dt + E \gamma(E ^*\bar{Y}_t)\,dt+ (\lambda -A)D_1\,d\tilde{W}_t + G(t,\bar{X}_t)\, dW_t &\\ \\
 - d{\bar{Y}}_t = A^T {\bar{Y}}_t \, dt   + l^0_x(t,\bar{X}_t) \, dt  -\bar{Z}_t \, dW_t  -\tilde{Z}_t \, d\tilde{W}_t , \qquad t \in [0,T] & \\ \\
\bar{X}_0 =x, \ \bar{Y}_T = - h_x(\bar{X}_T) ,\end{array}\right.
\end{equation*}
\end{remark}

\subsection{An auxiliary LQ control problem}\label{sez-LQ}
This section is devoted to the solution of the affine FBSDE. Let $b_0, h_0 \in  L^2_\P((0,T)\times \Omega; H)$ and $g_0 \in L^2 (\Omega, \F_T;H)$, consider: 
\begin{equation}\label{forwardback_linear}
\left\{\begin{array}{ll}  d\bar{X}_t = A \bar{X}_t  dt - [E + B][E+ B]^*\bar{Y}_t\,dt+ b_0(t)\,dt + (\lambda -A)D_1\,d\tilde{W}_t + G\, dW_t &\\ \\
 - d{\bar{Y}}_t = A^T {\bar{Y}}_t \, dt   +\bar{X}_t \, dt + h_0(t) \, dt  -\bar{Z}_t \, dW_t  -\tilde{Z}_t \, d\tilde{W}_t , \qquad t \in [0,T] & \\ \\
\bar{X}_0 =x, \, -\bar{Y}_T = \bar{X}_T + g_0 \end{array}\right.
\end{equation}
This system  is the Hamiltonian system corresponding to the control problem with state equation:
\begin{equation}\label{forward_lin}
\left
\{\begin{array}{ll}  dX_t = A X_t \, dt + [E+ B]u_t \,dt +  b_0(t) \, dt + (\lambda -A)D_1\,d\tilde{W}_t + G\, dW_t & t \in [0,T]\\
X_0 =x ,\end{array}
\right.
\end{equation} 
and cost functional 
\begin{equation}\label{funzionalecosto}
J(x,u)= \frac{1}{2} \E \int_0^T (|X_t + h_0(t)|^2 + |u_t|^2) \, dt +  \frac{1}{2} \E |X_T+g_0|^2
\end{equation}
to minimize over all $u \in \mathcal{U}$.
We will exploit this interpretation through the control problem in order to solve \eqref{forwardback_linear}, to this purpose we introduce the following Riccati equation:
\begin{equation}\label{Riccati}
\left\{\begin{array}{ll} 
\displaystyle - \frac{dP_t}{dt} = A^*P_t + P_t A - P_t (E+B) (E+B)^*P_t + I, & t \in [0,T] \\
\quad \, P_T= I
\end{array}\right.
\end{equation}
and the following backward equation, to cope with the affine terms:
 \begin{equation}\label{back_linear}
\left\{\begin{array}{ll} 
-dr_t = A^*r_t \,dt - P_t (E+B) (E+B)^*r_t\, dt+ P_t b_0(t) \,dt - h_0(t)\, dt - q_t \, W_t  -\tilde{q}_t \, d\tilde{W}_t, & t \in [0,T] \\
\quad r_T=I g_0
\end{array}\right.
\end{equation}
We denote, as in \cite{BenDel}, by $\Sigma(H)$ the space of self adjoint linear operators in $H$
and by $C_s([0,T]; \Sigma(H))$ the space of all
strongly continuous mappings from $[0,T]$ to $\Sigma(H)$, that is $P:[0,T]\rightarrow \Sigma(H)$
such that for every $h\in H$ $t\mapsto P_t h$ is continuous.

\noindent In the book \cite{BenDel} (part. IV, Chapter 2, Theorem 2.1),
it is proved that the first equation \eqref{Riccati} has a

solution in the space $ C_{s,\alpha}([0,T]; \Sigma(H))$,
the set of all $P \in  C_s([0,T]; \Sigma(H))$ such that:
\begin{itemize}
\item[(i)] $P(t)x \in D((-A^*)^{1-\alpha})$, for all $x \in H$, $t \in [0,T[$,
\item[(ii)] $(-A^*)^{1-\alpha}P \in C([0,T[; L(H))$,
\item[(iii)]$ \lim_{ t \to T} (T-t)^{1-\alpha} (-A^*)^{1-\alpha}P_t x=0$, for all $x \in H$. 
\end{itemize}
Moreover define
\begin{equation}
\| P \|_1= \sup_{ t \in [0,T[} \|  (T-t)^{1-\alpha} (- A^*)^{1-\alpha} P(t) \|
\end{equation}
$ C_{s,\alpha}([0,T]; \Sigma(H))$, endowed with the norm 
\begin{equation}
\| P \|_\alpha=  \| P\| + \| P\|_1
\end{equation}
is a Banach space.
We can now prove existence and uniqueness of a solution to \eqref{back_linear}, for semplicity we will denote the couple $(q,\tilde {q})$ as  $\hat{q}$ along with the comprehensive Wiener  process $\hat{W}_t := (W_t,\tilde{W}_t)$:
\begin{theorem}
Assume {\bf(A)} and  {\bf(B)} .  Then equation \eqref{back_linear} has a unique mild solution $(r_t, \hat{q}) \in L^2_\P (\Omega;C([0,T];H)) \times L^2_\P(\Omega \times [0,T]; L_2(\Xi \times K;H))$, moreover:
\begin{equation}\label{stima_back_linear}
\E\sup_{ t \in [0,T[} (T-t)^{2(1-\alpha)} |r_t| ^2 < \infty .
\end{equation} 
\end{theorem}

\Dim
We will prove existence and uniqueness by a fixed point technique.
Let us define a  map $\Gamma:  Y \to Y$, where $$Y := \left\lbrace   (r,\hat{q}) \in L^2_\P (\Omega;C([0,T];H)) \times L^2_\P(\Omega \times [0,T]; L_2(\Xi \times K;H)) : \E\sup_{ t \in [0,T[} (T-t)^{2(1-\alpha)} | r_t| ^2 < \infty\right\rbrace $$
such that $\Gamma((r', \hat q')) = (r, \hat q )$ is the {\em mild}  solution to:
\begin{equation}\label{mildptofix}\left\{
\begin{aligned}
r_t =  &e^{A^*(T-t)}g_0 - \int_t^T e^{A^* (s-t)}P_s (E+B) (E+B)^*r'_s\, ds +  \int_t^T e^{A^* (s-t)} P_s b_0(s) \,ds \\ & -  \int_t^T e^{A^* (s-t)}h_0(s)\, ds -   \int_t^T e^{A^* (s-t)}  \hat {q}_s\, d\hat{W}_s, & t \in [0,T] 
\end{aligned}
\right.
\end{equation}
We will prove that:
\begin{itemize}
\item [1)] $\Gamma ((r', \hat q')) \in Y$,
\item [2)] for any $\alpha < 1$ there exists  $\delta \in [0,T[$  that depends only on $\alpha$ and  constants appearing in ${\bf (A)}$ and ${\bf (B)}$ and $T$ such that
 \begin{equation}\label{ptofix}\| (r^1, \hat q^1)- (r^2, \hat q^2)\|_{Y_\delta} \leq  \alpha \| (r'^1, \hat q'^1)- (r'^2, \hat q'^2)  \|_{Y_\delta}\end{equation}
for some $\delta>0$ and we set
\begin{equation}
\begin{aligned}
Y_\delta:= \!  \left\lbrace   (r,\hat{q}) \in L^2_\P (\Omega;C([T-\delta,T];H)) \times L^2_\P(\Omega \times (T-\delta,T); L_2(\Xi \times K;H)): \right.\\\left.  \!\!\E\!\!\sup_{ t \in [T-\delta,T[} \!\!(T-t)^{2(1-\alpha)} | (\lambda - A^*)^{1-\alpha} r_t| ^2 < \infty \right\rbrace .
\end{aligned}
\end{equation}
\end{itemize}
The space $Y_\delta$ endowed with the norm:
$$
\| (r,\hat{q}) \|^2_{Y_\delta} := \E  \sup_{t \in [T-\delta,T]} |r_t| ^2 +  \E  \sup_{t \in [T-\delta,T[} (T-t)^{2(1-\alpha)} | (\lambda - A^*)^{1-\alpha} r_t| ^2 + \E\int_{T-\delta}^T |\hat{q}_t| ^2 \, dt
 $$
 is a Banach space. \\
  {\em Proof of statement 1)}.
 
 We introduce  the approximating problems for $ k > \lambda$:
\begin{equation}\label{Riccati_app}
\left\{\begin{array}{ll} 
\displaystyle - \frac{dP^k_t}{dt} = A^*P^k_t + P^k_t A - P^k_t (E^k+B) (E^k+B)^*P^k_t + I, & t \in [0,T] \\
\quad \, P^k_T= I
\end{array}\right.
\end{equation}
where $E^k := (\lambda - A)^{1-\alpha} kR(k,A) (\lambda - A)^{\alpha}D$, with $R(k,A):= (k -A)^{-1}$.

From \cite{BenDel} we know that  equation
 \eqref{Riccati_app} has a unique mild solution $P ^k \in C_{s,\alpha}([0,T]; \Sigma(H))$, for every $k$ and moreover the  following holds, see \cite{BenDel}(part IV, Chapter 2, lemma 2.1) :
\begin{equation}\label{Riccati_reg}
\left\{\begin{array}{ll} 
\lim_{k \to \infty} P^k(\cdot) x= P(\cdot) x \quad  \text{ in } C([0,T];H), \\
\lim_{k \to \infty} (T-\cdot)^{1-\alpha}(-A^*)^{1-\alpha}P^k(\cdot) x= (T-\cdot)^{1-\alpha}(-A^*)^{1-\alpha}P(\cdot) x \quad \text{ in }C([0,T];H).
\end{array}\right.
\end{equation}
Given $P^k$ we introduce also:
\begin{equation}\label{back_linear_app}
\left\{\begin{array}{ll} 
-dr^k_t = A^*r^k_t \,dt - P^k_t (E^k+B) (E^k+B)^*{r'}_t\, dt+ P^k_t b_0(t) \,dt - h_0(t)\, dt - \hat{q}^k_t d\, \hat{W}_t  & t \in [0,T] \\
\quad r^k_T=  g_0.
\end{array}\right.
\end{equation}
Existence and uniqueness of a mild solution for equation \eqref{back_linear} in $ L^2_\P (\Omega;C([0,T];H)) \times L^2_\P(\Omega \times (0,T); L_2(\Xi \times K;H))$ can be deduced by \cite{HuPeng}(prop. 2.1).
Now we can prove that 
\begin{align}\label{stima_L1}
& \displaystyle\E \int_0^T |\hat{q}^k_t|^2 \,dt\leq C \Big[E |g_0|^2 \\ \nonumber
&  + E \Big(\int_0^T |P^k_s (E^k+ B)(E^k+ B)^* {r'}_s|\, ds\Big)^2 +  E \int_0^T |P^k_tb_0(s)|^2\, ds+ 
E \int_0^T |h_0(s)|^2\, ds\Big]
\end{align}
The former estimate can be achieved  evaluating $d_t|r_t|^2$ and exploiting the fact that, being $A^*$ the generator of a contraction semigroup, $ \langle  A^* y,y\rangle \leq \omega |y|^2$, for any $ y \in D(A^*)$.  Since $r^k$ does not belong to $D(A^*)$, we  multiply $r^k$ by $nR(n,A)$, for $ n > \omega$ in order to perform the It\^o formula.  
Let us set $r^{n,k}_t = nR(n,A) r^k_t$ and $\hat{q}^{n,k} =nR(n,A) \hat{q}^{k}_t$ ,  hence:
\begin{equation}\label{back_linear_app_bis}
\left\{\begin{aligned}
-dr^{n,k}_t = & A^*r^{n,k}_t \,dt - nR(n,A)P^k_t (E^k+B) (E^k+B)^*r'_t\, dt+ nR(n,A)P^k_t b_0(t) \,dt\\ &  - nR(n,A)h_0(t)\, dt - \hat{q}^{n,k}_t d\, \hat{W}_t , \quad t \in [0,T] \\
\quad r^{n,k}_T= & nR(n,A) g_0.
\end{aligned}\right.
\end{equation}
Now we can evaluate $d_t |r^{n,k}_t|^2$:
\begin{equation}
\begin{aligned}
d_t |r^{n,k}|^2= & 2 \langle  A^* r^{n,k}_t, r^{n,k}_t\rangle \, dt  -  2 \langle f^{n,k}_t,  r^{n,k}_t\rangle  \, dt - 2 \langle \hat{q}^{n,k}_t,   r^{n,k}_t\rangle \, d\hat{W}_t - |q^{n,k}_t|^2 \, dt
\end{aligned}
\end{equation}
where 
\[  f^{n,k}_t=  nR(n,A) P^k_t (B+E^k)(B+E^k)^*{r'}_t + P^k_t b_0(t) + h_0(t)
\]
Now similarily to  \cite{Gua1} (prop. 3.4), see also \cite{Magnifici5}(lemma 3.1), we get:
\begin{equation} \label{stimaL2nk}
\E \int_0^T |\hat{q}^{n,k}_t|^2 \, dt  \leq C \Big [  \E \sup_{ t \in [0,T]} |r ^{n,k}_t|^2 + \E \Big (  \int_0^T |f ^{n,k}_t| \, dt \Big )^2\Big ]
 \end{equation}
 where the constant $C$ depends on constants appearing in ${\bf (A)}$ and ${\bf (B)}$ and $T$.
 Letting $n$ tend to $\infty$ we obtain estimate \eqref{stima_L1}.
 Now bearing in mind that $\sup_{t \in [0,T]}|P^k_t| \leq M$ independent of $k$, thanks to  \eqref{Riccati_reg}  and Banach-Steinhaus theorem, we obtain that:
 \begin{equation}\label{stima_L1_bis}
 \begin{aligned}
 \displaystyle\E \int_{0}^T |\hat{q}^k_t|^2 \,dt  \leq &C \Big[\E \Big(\sup_{ s \in [0,T]} (T-s)^{(2-2\alpha)}[|(\lambda - A^*)^{1-\alpha} {r}'_s|^2 + |r'_s|^2]\, ds\int_0^T  s^{\alpha -1} (T-s)^{2\alpha -2}\, ds\Big)^2 \\ &
 +\E |g_0|^2 +   \E \int_0^T |b_0(s)|^2\, ds+ 
\E \int_0^T |h_0(s)|^2\, ds\Big]
\end{aligned}
\end{equation}
Let us consider $k,m > \omega$:
\begin{align*}
r^k_t-r^m_t &= - \int_0^T e^{A^*(s-t)} [P^k_s (B+E^k)(B+E^k)^*- P^m_s (B+E^m)(B+E^m)^*]{r'}_s\, ds \\
&-  \int_0^T e^{A^*(s-t)}  (\hat{q}^k_s - \hat{q}^m_s) \, d\hat{W}_s.
\end{align*}
We have that:
 \begin{align*}
&r^k_t-r^m_t =  \Et (r^k_t-r^m_t) \\
&= - \Et \int_0^T e^{A^*(s-t)} [P^k_s (B+E^k)(B+E^k)^*- P^m_s (B+E^m)(B+E^m)^*]{r'}_s\, ds,
\end{align*}
and, since $| (\lambda - A^*)^{1-\alpha}e ^{s A^*}|_{L(H)} \leq c\,s^{1-\alpha}$:
\begin{equation*}
\begin{aligned}
&\sup_{ t \in [0,T[}|(\lambda - A^*)^{1-\alpha}(r^k_t-r^m_t)|^2
\\& \leq c\, \Et \int_0^T s^{\alpha -1}|[P^k_s (B+E^k)(B+E^k)^*-P(B+E)(B+E)^*]r'_s| \, ds \\  & +c\, \Et\int_0^T s^{\alpha -1}|[P^m_s (B+E^m)(B+E^m)^*-P(B+E)(B+E)^*]r'_s| \, ds
\end{aligned}
\end{equation*}
Hence, taking into account that:
\begin{equation}
\E \sup_{s \in [0,T[ }(T-s)^{2(1-\alpha)}|(\lambda - A^*)^{1-\alpha}r'_s|^2 < \infty ,
\end{equation}
by dominated convergence we end up with
\begin{equation}
\lim_{k,m \to + \infty} \E\sup_{ t \in [0,T[} (T-t)^{2(1-\alpha)}|(\lambda - A^*)^{1-\alpha}(r^k_t-r^m_t)|^2 =0
\end{equation}
Similarily we have that  
\begin{equation}
\lim_{k,m \to + \infty} \E|(r^k_t-r^m_t)|^2 =0
\end{equation}
Moreover from former calculations we have that
 \begin{equation}\label{stima_L1_cauchy}
 \begin{aligned}
 \displaystyle\E \int_{0}^T |\hat{q}^k_t -\hat{q}^m_t  |^2 \,dt  \leq &C \Big[ \E\Big ( \int_0^T  |[P^k_s (B+E^k)(B+E^k)^*-P(B+E)(B+E)^*]r'_s | \, ds \Big)^2 \\ &
  +\E\Big (\int_0^T |[P^m_s (B+E^m)(B+E^m)^*-P(B+E)(B+E)^*]r'_s| \, ds \Big)^2 \Big]
\end{aligned}
\end{equation}
Thus, the limit processes $r$ and $\hat{q}$  solve equation \eqref{mildptofix} and we have the desired regularity.

Now we have to prove \eqref{ptofix}. Following previous procedures, 
we have:
\begin{align*}
&\E \sup_{t \in [T-\delta,T[} (T-t)^{2(1-\alpha)}| (\lambda - A^*)^{1-\alpha}| (r^1_t - r^2_t)| ^2 \\
&\leq  M \delta^{2(2\alpha -1)}
\E \sup_{t \in [T-\delta,T[} (T-t)^{2(1-\alpha)}| (\lambda - A^*)^{1-\alpha}| (r'^1_t - r'^2_t)| ^2
\end{align*}
\begin{equation*}
\E \sup_{t \in [T-\delta,T]} | r^1_t - r^2_t| ^2 \leq  M \delta^{2(2\alpha -1)}
\E \sup_{t \in [T-\delta,T[} |(\lambda - A^*)^{1-\alpha} (r'^1_t - r'^2_t)| ^2
\end{equation*}
and 
 \begin{equation*}
 \begin{aligned}
 \displaystyle\E \int_{T-\delta}^T [|\hat{q}^1_t -\hat{q}^2_t  |^2 \,dt  &\leq M \Big[ \E\Big ( \int_0^T  |P_t (B+E)(B+E)^*(r'^1_t-r'^2_t) | \, dt \Big)^2  \Big] \\ &
 \leq M \delta^{2(2\alpha -1)} \E \sup_{t \in [T-\delta,T[} |(\lambda - A^*)^{1-\alpha} (r'^1_t - r'^2_t)| ^2
\end{aligned}
\end{equation*}
where the constant $M$ depends on $\alpha$ and constants appearing in hypotheses ${\bf (A)}$ and ${\bf (B)}$
and $M\delta^{2\alpha-1}<1$ if $\delta$ is sufficiently small. Therefore one can repeat the procedure in $[T-2\delta, T-\delta]$ and so on in order to cover, in a finite number of steps, the whole interval $[0,T]$. \finedim

It remains to show that if we define $\bar{Y}_t=P_t \bar{X}_t+r_t$, then $\bar{Y} $ is a solution
to the BSDE in the FBSDE (\ref{forwardback_linear}).
\begin{proposition}\label{prop-fbsde-lin}
Let assumptions $\textbf{(A)}$ hold true an let $b_0\,,h_0\in L^2_\P (\Omega\times[0,T];H),
\,g_0\in L^2(\Omega;H)$. Then the FBSDE (\ref{forwardback_linear}) admits a
unique mild solution $(\bar X, \bar Y, (\bar Z, \tilde Z))
\in L^2_\P (\Omega;C([0,T];H)) \times L^2_\P (\Omega;C([0,T];H)) 
\times L^2_\P(\Omega \times [0,T]; L_2(\Xi \times K;H))$ satisfying moreover $\E\sup_{t\in[0,T]}(T-t)^{2(1-\alpha)}\Vert (E+B)\bar{Y}_t\Vert^2<+\infty$.
\end{proposition}
\Dim Let us denote by $P^k$ the solution of the Riccati equation (\ref{Riccati_app}) and,
for $j>\omega$, by $A_j:=jR(j,A)$
the Yosida approximants of $A$. We denote by $P^{j,k}$ the solution of the Riccati equation
(\ref{Riccati_app}) with $A_j$ in the place of $A$:
\begin{equation}\label{Riccati_app_bis}
\left\{\begin{array}{ll} 
\displaystyle - \frac{dP^{j,k}_t}{dt} = A_j^*P^{j,k}_t + P^{j,k}_t A - P^{j,k}_t (E^k+B)
 (E^k+B)^*P^{j,k}_t + I, & t \in [0,T] \\
\quad \, P^{j,k}_T= I
\end{array}\right.
\end{equation}
By $(r^k,(g^k,\tilde{g}^k))$ and by $(r^{n,k},(g^{n,k},\tilde{g}^{n,k}))$
we denote respectively the solution of the BSDEs (\ref{back_linear_app})
and (\ref{back_linear_app_bis}).
Moreover we denote by $\bar{X}$ and $\bar{X}^k$ respectively the solution of
\begin{equation}\label{forward_linear_P}
\left\{\begin{array}{ll}  d\bar{X}_t = A \bar{X}_t  dt - [E + B][E+ B]^*
(P_t \bar{X}_t+r_t)\,dt+ b_0(t)\,dt + (\lambda -A)D_1\,d\tilde{W}_t + G\, dW_t &\\ \\
\bar{X}^{k}_0=x, \end{array}\right.
\end{equation}
and of
\begin{equation}\label{forward_linear_app}
\left\{\begin{array}{ll}  d\bar{X}^{k}_t = A \bar{X}^{k}_t \, dt - [E^k + B][E^k+ B]^*
(P^{k}_t \bar{X}^{k}_t+r^{k}_t)\,dt+ b_0(t)\,dt + (\lambda -A)D_1\,d\tilde{W}_t + G\, dW_t &\\ \\
\bar{X}^{k}_0=x. \end{array}\right.
\end{equation}
We also set $\bar{X}^{n,k}=nR(n,A)\bar{X}^{n,k}$ which is solution of
\begin{equation}\label{forward_linear_app_bis}
\left\{\begin{array}{ll}  d\bar{X}^{n,k}_t = A_n \bar{X}^{k}_t \, dt - nR(n,A)[E^k + B][E^k+ B]^*
(P^{k}_t \bar{X}^{k}_t+r^{k}_t)\,dt+ b_0(t)\,dt &\\ \\
\qquad\qquad +nR(n,A)(\lambda -A)D_1\,d\tilde{W}_t +nR(n,A) G(t,\bar{X}_t)\, dW_t &\\ \\
\bar{X}^{k}_0=nR(n,A)x. \end{array}\right.
\end{equation}
By applying It\^o formula to $P^{j,k}_t\,\bar{X}^{n,k}_t+r^{n,k}_t$
we get
\begin{align*}
& d (P^{j,k}_t\,\bar{X}^{n,k}_t+r^{n,k}_t)=\left(-A_j^*P^{j,k}_t\bar{X}^{n,k}_t
-P^{j,k}_t A\bar{X}^{n,k}_t +P^{j,k}_t[E^k + B][E^k+ B]^*P^{j,k}_t\,\bar{X}^{n,k}_t+\bar{X}^{n,k}_t\right) \,dt\\ &+P^{j,k}_tA_j\bar{X}^{n,k}_t\,dt
-P^{j,k}_t nR(n,A)[E^k + B][E^k+ B]^*\left(P^{k}_t\,\bar{X}^{k}_t+r^{k}_t\right)\,dt
+P^{j,k}_t nR(n,A)b_0(t)\,dt\\
& +P^{j,k}_t nR(n,A)(\lambda -A)D_1\,d\tilde{W}_t +P^{j,k}_t nR(n,A) G\, dW_t -A^*r^{n,k}_t\,dt
+nR(n,A)h_0(t)\,dt\\
&
+nR(n,A)P^k_t[E^k + B][E^k+ B]^*r^{n,k}_t\,dt-nR(n,A)P^k_t b_0(t)\,dt +q^{n,k}_t\,dW_t
+\tilde{q}^{n,k}_t\,d\tilde{W}_t.
\end{align*}
So in mild form we get
\begin{align*}
 & P^{j,k}_t\,\bar{X}^{n,k}_t+r^{n,k}_t=e^{(T-t)A^*_j}[\bar{X}^{n,k}_T+nR(n,A)g_0]
 +\int_t^T e^{(s-t)A^*_j}\left(A^*r^{n,k}_s-A^*_jr^{n,k}_s\right)\,ds\\
&+\int_t^T e^{(s-t)A^*_j} \left(nR(n,A)P_s b_0(s)-P^{j,k}_s nR(n,A)b_0(s)\right)\,ds
 +\int_t^T e^{(s-t)A^*_j}P^{j,k}_s\left(A_j\bar{X}^{n,k}_s-A\bar{X}^{n,k}_s\right)\,ds\\
 &- \int_t^T e^{(s-t)A^*_j}\left[P^{j,k}_s[E^k + B][E^k+ B]^*(P^{j,k}_s\,\bar{X}^{n,k}_s+r^{n,k}_s)
 +\bar{X}^{n,k}_s\right] \,ds\\
&+\int_t^T e^{(s-t)A^*_j}P^{j,k}_s nR(n,A)[E^k + B][E^k+ B]^*(P^{k}_s\,\bar{X}^{k}_s+r^k_s) \,ds-\int_t^T e^{(s-t)A^*_j}nR(n,A)h_0(s)\,ds\\
&\int_t^T e^{(s-t)A^*_j}\left(P^{j,k}_s nR(n,A)(\lambda -A)D_1-\tilde{q}^{n,k}_s\right)\,d\tilde{W}_t
 +\int_t^T e^{(s-t)A^*_j}\left(P^{j,k}_s nR(n,A) G-q^{n,k}_t\right)\, dW_t .
\end{align*}
We start by letting $j\rightarrow \infty$. It follows by assumption $\textbf{(A.1)}$ that
$\Vert e^{tA_j} \Vert\leq e^{\omega t} $. Keeping this in mind, and since $r^{n,k},\bar{X}^{n,k}
\in \cald(A)$ and moreover since $P^{j,k}$ is uniformly bounded in $j$, we get that the integrals
$\displaystyle\int_t^T e^{(s-t)A^*_j}\left(A^*r^{n,k}_s-A^*_jr^{n,k}_s\right)\,ds$
and $\displaystyle\int_t^T e^{(s-t)A^*_j}P^{j,k}_s\left(A_j\bar{X}^{n,k}_s-A\bar{X}^{n,k}_s\right)\,ds$
converge to $0$ as $j\rightarrow \infty$.

\noindent With similar considerations, by adding and subtracting $e^{(s-t)A^*_j}P^{j,k}_s[E^k + B][E^k+ B]^*P^{k}_s\,\bar{X}^{n,k}_s$ and $e^{(s-t)A^*_j}P^{k}_s[E^k + B][E^k+ B]^*P^{k}_s\,\bar{X}^{n,k}_s$
we get that the integral $\displaystyle\int_t^T e^{(s-t)A^*_j}P^{j,k}_s[E^k + B][E^k+ B]^*P^{j,k}_s\,\bar{X}^{n,k}_s\,ds$  converges to $\displaystyle\int_t^T e^{(s-t)A^*}P^{k}_s[E^k + B][E^k+ B]^*P^{k}_s\,\bar{X}^{n,k}_s\,ds$ as $j\rightarrow \infty$.

\noindent In an analogous and simpler way we also get that $\displaystyle 
\int_t^T e^{(s-t)A^*_j}P^{j,k}_s nR(n,A)[E^k + B][E^k+ B]^*(P^{k}_s\,\bar{X}^{k}_s+r^k_s) \,ds$
converges to $\displaystyle \int_t^T e^{(s-t)A^*}P^{k}_s 
nR(n,A)[E^k + B][E^k+ B]^*(P^{k}_s\,\bar{X}^{k}_s+r^k_s) \,ds$.

\noindent By adding and subtracting $P^{j,k}_s nR(n,A)b_0(s)$ it is possible to see that $$\displaystyle\int_t^T e^{(s-t)A^*_j} \left(nR(n,A)P_s b_0(s)-P^{j,k}_s nR(n,A)b_0(s)\right)\,ds$$ converges to 
$\displaystyle\int_t^T e^{(s-t)A^*_j} \left(nR(n,A)P_s b_0(s)-P^{k}_s nR(n,A)b_0(s)\right)\,ds$ as $j\rightarrow \infty$; it is immediate to see that $\displaystyle\int_t^T e^{(s-t)A^*_j}nR(n,A)h_0(s)\,ds
\rightarrow\displaystyle\int_t^T e^{(s-t)A^*}nR(n,A)h_0(s)\,ds$ as $j\rightarrow \infty$..

\noindent For what concerns the stochastic integrals, we notice that the integrands
are square integrable with respect to $s$, uniformly with respect to $j$.

\noindent So, letting $j\rightarrow \infty$, we get
\begin{align*}
 & P^{k}_t\,\bar{X}^{n,k}_t+r^{n,k}_t=e^{(T-t)A^*}[\bar{X}^{n,k}_T+nR(n,A)g_0]\\
&+\int_t^T e^{(s-t)A^*} \left(nR(n,A)P^k_s b_0(s)-P^{k}_s nR(n,A)b_0(s)\right)\,ds\\
&- \int_t^T e^{(s-t)A^*}\left[P^{k}_s[E^k + B][E^k+ B]^*(P^{k}_s\,\bar{X}^{n,k}_s+r^{n,k}_s)
 +\bar{X}^{n,k}_s\right] \,ds-\int_t^T e^{(s-t)A^*}nR(n,A)h_0(s)\,ds\\
&+\int_t^T e^{(s-t)A^*}P^{k}_s nR(n,A)[E^k + B][E^k+ B]^*(P^{k}_s\,\bar{X}^{k}_s+r^k_s) \,ds\\
&\int_t^T e^{(s-t)A^*}\left(P^{k}_s nR(n,A)(\lambda -A)D_1-\tilde{q}^{n,k}_s\right)\,d\tilde{W}_t
 +\int_t^T e^{(s-t)A^*}\left(P^{k}_s nR(n,A) G-q^{n,k}_s\right)\, dW_s .
\end{align*}
Moreover,
\begin{align*}
 &\Et P^{k}_t\,\bar{X}^{n,k}_t+r^{n,k}_t= P^{k}_t\,\bar{X}^{n,k}_t+r^{n,k}_t\\
&=\Et e^{(T-t)A^*}[\bar{X}^{n,k}_T+nR(n,A)g_0]
+\Et\int_t^T e^{(s-t)A^*} \left(nR(n,A)P^k_s b_0(s)-P^{k}_s nR(n,A)b_0(s)\right)\,ds\\
&- \Et\int_t^T e^{(s-t)A^*}\left[P^{k}_s[E^k + B][E^k+ B]^*(P^{k}_s\,\bar{X}^{n,k}_s+r^{n,k}_s)
 +\bar{X}^{n,k}_s\right] \,ds-\int_t^T e^{(s-t)A^*}nR(n,A)h_0(s)\,ds\\
&+\Et\int_t^T e^{(s-t)A^*}P^{k}_s nR(n,A)[E^k + B][E^k+ B]^*(P^{k}_s\,\bar{X}^{k}_s+r^k_s) \,ds.\\
\end{align*}
As $n\rightarrow \infty$, taking into account, where necessary, that $\Vert P^k
\Vert_1$ is bounded uniformly with respect to $k$, we get
\begin{align*}
 &\Et P^{k}_t\,\bar{X}^{k}_t+r^{k}_t=\Et e^{(T-t)A^*}[\bar{X}^{k}_T+g_0]\\
&+\Et\int_t^T e^{(s-t)A^*} \left(P^k_s b_0(s)-P^{k}_s b_0(s)\right)\,ds\\
&- \Et\int_t^T e^{(s-t)A^*}\left[P^{k}_s[E^k + B][E^k+ B]^*(P^{k}_s\,\bar{X}^{k}_s++r^{k}_s)
 +\bar{X}^{k}_s\right] \,ds+\Et\int_t^T e^{(s-t)A^*}(\bar{X}^{k}_s-h_0(s))\,ds\\
&+\Et\int_t^T e^{(s-t)A^*}P^{k}_s [E^k + B][E^k+ B]^*(P^{k}_s\,\bar{X}^{k}_s+r^k_s) \,ds\\
&=\Et e^{(T-t)A^*}[\bar{X}^{k}_T+g_0]+\Et\int_t^T e^{(s-t)A^*}(\bar{X}^{k}_s-h_0(s))\,ds.\\
\end{align*}

Now notice that
\begin{equation*}
 \bar{X}^{k}_t-\bar{X}_t=\int_0^t e^{(t-s)A}\left(  [E^k + B][E^k+ B]^*
(P^{k}_s \bar{X}^{k}_s+r^{k}_S)-[E + B][E+ B]^*
(P_s \bar{X}_s+r_s)\right)\,ds.
\end{equation*}
By the convergence of $P^k$ to $P$, see \cite{BenDel}, chapter IV, section 2, lemma 2.1
and theorem 2.1, by adding and subtracting suitable terms and in virtue of Gronwall lemma,
we get that $ \bar{X}^{k}\rightarrow\bar{X} $ in $L^2_\P(\Omega,C([0,T],H))$.
By adding and subtracting $P^{k}_t\,\bar{X}_t$ we also get that
$P^{k}_t\,\bar{X}^k_t\rightarrow P_t\,\bar{X}_t$ in $L^2_\P(\Omega,C([0,T],H))$,
since $\sup_{t \in [0,T]}|P^k_t| \leq M$ independent of $k$, thanks to  \eqref{Riccati_reg}  and Banach-Steinhaus theorem. With similar arguments we also get that $(T-t)^\alpha(\lambda-A^*)^{1-\alpha}P^k_t\,\bar{X}^k_t$
converges to $(T-t)^\alpha(\lambda-A^*)^{1-\alpha}P_t\,\bar{X}_t$
 in $L^2_\P(\Omega,C([0,T],H))$.
With similar and simpler arguments we finally get
\begin{equation*}
  P_t\,\bar{X}_t+r_t=\Et e^{(T-t)A^*}[\bar{X}_T+g_0]\\
- \Et\int_t^T e^{(s-t)A^*}\bar{X}_s \,ds-\Et\int_t^T e^{(s-t)A^*}h_0(s)\,ds\\
\end{equation*}
Arguing as in \cite{HuPeng} and as in the proof of lemma \ref{lemma_bridge}, by the extended martingale representation theorem, (see also
\cite{HuPeng1} and \cite{Yor})), for every $s\in[0,T]$ there exists $(K(s,\cdot),\tilde{K}(s,\cdot),
\in L^2_\P(\Omega\times[0,T],L_2(\Xi,H))\times L_\P^2(\Omega\times[0,T],L_2(K,H))$
such that $\forall\, 0\leq t\leq s\leq T$
\begin{equation*}
 \Et \int_t^T e^{(s-t)A^*}(\bar{X}_s-h_0(s))\,ds=\E \int_t^T e^{(s-t)A^*}(\bar{X}_s
h_0(s))\,ds+
\int_0^t K(s,\theta )\,dW_\theta+
\int_0^t \tilde{K}(s,\theta )\,d\tilde{W}_\theta
\end{equation*}
Note that $\forall \,\theta\geq s$, $K(s,\theta)=\tilde{K}(s,\theta)=0$ a.e. and
\begin{equation}
\label{stimaKbis}
 \E \int_0^T\int_0^s\left[\vert K^{n,k}(s,\theta)\vert^2
+\vert \tilde{K}^{n,k}(s,\theta)\vert^2 \right]\, d\theta \,ds
\leq 4\E \int_0^T \vert (\bar{X}^{n,k}_s- nR(n,A)h_0(s))\vert^2 \,ds.
\end{equation}
Moreover, there exists $(L,\tilde{L}) \in 
L^2_\P(\Omega\times[0,T],L_2(\Xi,H))\times L^2_\P(\Omega\times[0,T],L_2(K,H))$ such that 
\[
\Et [\bar{X}_T+g_0]=\E[\bar{X}_T+g_0]+\int_0^t L(\theta )\,dW\theta+
\int_0^t \tilde{L}(\theta )\,d\tilde{W}_\theta.
\]
We deduce that by setting, for almost all $s\in[0,T]$,
\begin{equation*}
Z_s=\int_s^Te^{(\alpha-s)A^T}K_s(\alpha,s)\,d\alpha,\quad
\tilde{Z}_s=\int_s^Te^{(\alpha-s)A^T} \tilde{K}_s(\alpha,s)\,d\alpha.
\end{equation*}
By the definition of $(K,\tilde{K}$ and
by estimates (\ref{stimaKbis}), it follows that
$(\bar{Z}_s, \tilde{Z}_s)\in
L^2_\P(\Omega \times [0,T]; L_2(\Xi \times K;H))$. 
Moreover $(X,Y,(Z,\tilde{Z}))$ are a solution to FBSDE (\ref{forwardback_linear})
\finedim
\subsection{Existence and uniqueness of the mild solution of the FBSDE}
In this section we prove theorem \ref{FBSDE}, by using th results in lemma \ref{lemma_bridge}
and in section \ref{sez-LQ}.

{\bf{Proof of Theorem \ref{FBSDE}}} \emph{Existence.} We follow the proof of Theorem 3.1, existence part, in \cite{HuPeng3},
with suitable changes due to the different framework.
For $\alpha \in[0,1]$ consider the FBSDE
\begin{equation}\label{forwardback_bridge_fin}
\left\{\begin{array}{ll}  d\bar{X}_t = A \bar{X}_t , dt+b^\alpha(\bar{Y}_t)\,dt+ b_0(t)\,dt + (\lambda -A)D_1\,d\tilde{W}_t + G\, dW_t &\\ \\
 - d{\bar{Y}}_t = A^* \bar{Y}_t \, dt   +h^\alpha(\bar{X}_t) \, dt + h_0(t) \, dt
  -\bar{Z}^{n}_t \, dW_t  -\tilde{Z}^{n}_t \, d\tilde{W}_t , \qquad t \in [0,T] & \\ \\
\bar{X}_0 =x, \, -\bar{Y}_T = g^\alpha(\bar{X}_t)+g_0\end{array}\right.
\end{equation}
For $\alpha=0$ the FBSDE
(\ref{forwardback_bridge_fin}) admits a mild solution: by
section \ref{sez-LQ} we know that FBSDE (\ref{forward_lin}) admits a mild solution,
and for $\alpha=0$, FBSDE
(\ref{forwardback_bridge_fin}) coincides with FBSDE (\ref{forward_lin}).
By lemma \ref{lemma_bridge} there exists $\delta_0$ such that
for all $\alpha\in [0,\delta_0]$ the FBSDE (\ref{forwardback_bridge_fin})
admits a mild solution with the required regularity. Then, by arbitrary
choice of $g_0,\,h_0 $ and $g_0$ we can solve (\ref{forwardback_bridge_fin})
for $\alpha\in [\delta_0,2\delta_0],\,[2\delta_0,3\delta_0],...$: notice that
$\delta_0$ does not depend on $\alpha$. We arrive at solving (\ref{forwardback_bridge_fin}) for $\alpha=1$, and again by the arbitrary choice we can make of  $g_0,\,h_0 $ and $g_0$
we have proved the existence of an adapted soltuion of (\ref{forwardback_ottimo})
$(\bar{X}, \bar{Y}, (\bar{Z}, \tilde{Z}))$ with the required regularity.

\noindent \emph{Uniqueness.} In order to prove uniqueness we follow
\cite{HuPeng3}, theorem 3.1 uniqueness part, and the proof of lemma
\ref{lemma_bridge} in the present paper.
Let, for $i=1,2$, 
$(\bar{X}^{i}, \bar{Y}^{i}, (\bar{Z}^{i}, \tilde{Z}^{i}))$ be
 two solutions of (\ref{forwardback_ottimo}). In order to apply It\^o formula, we have to approximate these solutions with elements in the
domain of $A$, namely we set  $(\bar{X}^{n,i}, \bar{Y}^{n,i}, (\bar{Z}^{n,i}, \tilde{Z}^{n,i}))
=(nR(n,A)\bar{X}^{i}, nR(n,A)\bar{Y}^{i}, (nR(n,A)\bar{Z}^{i}, nR(n,A)\tilde{Z}^{i}))$,
$i=1,2$, and as in lemma
\ref{lemma_bridge}, e also denote $E_n+B_n:=nR(n,A)(E+B)$.
By applying It\^o formula to
$\langle \bar{X}^{n,1}_t-\bar{X}^{n,2}_t,\bar{Y}^{n,1}_t-\bar{Y}^{n,2}_t \rangle$,
and then integrating over $[0,T]$ and taking
expectation we get
\begin{align}\label{uniq-ito-n}
 &-\E\< \bar{X}^{n,1}_T-\bar{X}^{n,2}_T,nR(n,A)(h_x(\bar{X}^{1}_T)-h_x(\bar{X}^{2}_T))\>
 \\ \nonumber
&=\E\int_0^T \< [E_n+B_n](\gamma([E+B]^* \bar{Y}_t^{1})-
\gamma([E+B]^* \bar{Y}_t^{2})),\bar{Y}_t^{n,1}-\bar{Y}_t^{n,2}\> \,dt\\ \nonumber
&-\E\int_0^T\< nR(n,A)(l^0_x(t,\bar{X}^{1}_t)-l^0_x(t,\bar{X}^{2}_t)),
\bar{X}^{n,1}_t-\bar{X}^{n,2}_t \> \,dt\\ \nonumber
\end{align}
Next we want to let $n\rightarrow +\infty$ in the (\ref{uniq-ito-n}): arguing as in lemma \ref{lemma_bridge},
we deudce that $\bar{X}^{n,i}\rightarrow \bar{X}^{n,i}$ in $L^2_\P(\Omega;C( [0,T], H))$ for $i=1,2$,
and that $\bar{Y}^{n,i}\rightarrow \bar{Y}^{n,i}$ in $L^2_\P(\Omega;C( [0,T], H))$ for $i=1,2$ and moreover
$\E\sup_{t\in[0,T]}(T-t)^{2(1-\alpha)}\vert [E+B]^*(\bar{Y}^{n,i}-\bar{Y}^{i})\vert^2\rightarrow 0$, for $i=1,2$.
We also have $\E\sup_{t\in[0,T]}(T-t)^{2(1-\alpha)}\vert [E_n+B_n]^*\bar{Y}^{n,i}-[E+B]^*\bar{Y}^{i})\vert^2\rightarrow 0$
for $i=1,2$. So letting $n\rightarrow \infty$ in (\ref{uniq-ito-n}) we get
\begin{align*}
 &-\E\< \bar{X}^{1}_T-\bar{X}^{2}_T,h_x(\bar{X}^{1}_T)-h_x(\bar{X}^{2}_T)\>
 \\ \nonumber
&=\E\int_0^T \< [E+B](\gamma([E+B]^* \bar{Y}_t^{1})-
\gamma([E+B]^* \bar{Y}_t^{2})),\bar{Y}_t^{1}-\bar{Y}_t^{2}\> \,dt\\ \nonumber
&-\E\int_0^T \< l^0_x(t,\bar{X}^{1}_t)-l^0_x(t,\bar{X}^{2}_t),
\bar{X}^{1}_t-\bar{X}^{1}_t \> \,dt.\\ \nonumber
\end{align*}
So, by assumptions ${\bf(B)}$ we get
\[
 \E\vert \bar{X}^{1}_T-\bar{X}^{2}_T \vert^2
+\E\int_0^T\vert[E+B]^*( \bar{Y}_t^{1}-\bar{Y}_t^{2})\vert^2\, dt +\E\int_0^T\vert \bar{Y}^{1}_t-\bar{Y}_t^{2}\vert^2\, dt\leq 0
\]
and so the uniqueness follows.
\finedim

\end{document}